\documentclass[aos,MSNbibl,seceqn,nameyear,dvips]{arximspdf}
\usepackage{graphicx}
\doi{10.1214/13-AOS1122} 
\volume{41}
\issue{3}
\pubyear{2013}
\firstpage{1642}
\lastpage{1668}

\makeatletter
\newcommand{\rrVert}{\Vert}
\newcommand{\rrvert}{\vert}
\newcommand{\llVert}{\Vert}
\newcommand{\llvert}{\vert}

\newtheorem{theorem}{Theorem}[section]
\newtheorem{lemma}{Lemma}[section]

\newproclaim{example}{Example}[section]
\newproclaim{condition}{Assumption}[section]
\newproclaim{definition}{Definition}[section]

\makeatother

\begin{document}
\begin{frontmatter}

\title{Regressions with Berkson errors in covariates---A~nonparametric approach}
\runtitle{Berkson errors}

\begin{aug}
\author[A]{\fnms{Susanne M.} \snm{Schennach}\corref{}\thanksref{t1}\ead[label=e1]{smschenn@alum.mit.edu}}
\runauthor{S. M. Schennach}
\affiliation{Brown University}
\address[A]{Department of Economics\\
Box B\\
Brown University\\
Providence, Rhode Island 02912\\
USA\\
\printead{e1}} 
\end{aug}

\thankstext{t1}{Supported in part by NSF Grants SES-0752699 and SES-1156347, and through TeraGrid computer
resources provided by the University of Texas under Grant SES-070003.}

\received{\smonth{8} \syear{2012}}
\revised{\smonth{4} \syear{2013}}

%
\begin{abstract}
This paper establishes that so-called instrumental variables enable the
identification and the estimation of a fully nonparametric regression
model with Berkson-type measurement error in the regressors. An
estimator is proposed and proven to be consistent. Its practical
performance and feasibility are investigated via Monte Carlo
simulations as well as through an epidemiological application
investigating the effect of particulate air pollution on respiratory
health. These examples illustrate that Berkson errors can clearly not
be neglected in nonlinear regression models and that the proposed
method represents an effective remedy.
\end{abstract}

%

\end{frontmatter}

\section{Introduction}

Many statistical data sets involve covariates $X$ that are
error-contaminated versions of their true unobserved counterpart
$X^{\ast}$. However, the measurement error often does not fit the
classical error
structure $X=X^{\ast}+\Delta X$ with $\Delta X$ independent from
$X^{\ast}$. A common occurrence is, in fact, the opposite situation,
in which $%
X^{\ast}=X+\Delta X^{\ast}$ with $\Delta X^{\ast}$ independent from $X$,
a situation often referred to as Berkson measurement error [\citet
{berskonme}, \citet{wangnlmbme}, \citet{CarrollME}]. A
typical example is
an epidemiological study in which an individual's true exposure
$X^{\ast}$
to some contaminant is not observed, but instead, what is available is the
average concentration $X$ of this contaminant in the region where the
individual lives. The individual-specific $X^{\ast}$ randomly fluctuate
around the region average $X$, resulting in Berkson errors.

Existing approaches to handle data with Berkson measurement error
[e.g., \citet
{delaiglenpberkson}, \citet{carrollberksonmix}] unfortunately
require the
distribution of the measurement error to be known, or to be estimated via
validation data, which can be costly, difficult or impossible to collect.
(In classical measurement error problems, the distribution of the error can
be identified from repeated measurements via a Kotlarski-type equality
[\citet
{schennachnlme}, \citet{livuong1998}]. However, such results do
not yet
exist for Berkson-type measurement error.) A popular approach to relax the
assumption of a fully known distribution of the measurement error is to
allow for some adjustable parameters in the distributions of the variables
and their relationships, and solve for the parameter values that best
reproduce various conditional moments of the observed variables, under the
assumption that this solution is unique. This approach has been used, in
particular, for polynomial specifications [\citet{huwangpoly}]
and, more
recently, for a very wide range of parametric models [Wang
(\citeyear{wangnlmbme,wangunified})].

The present paper goes beyond this and provides a formal identification
result and a general nonparametric regression method that is consistent in
the presence of Berkson errors, without requiring the distribution of the
measurement error to be known a priori. Instead, the method relies on the
availability of a so-called instrumental variable [e.g., see Chapter 6 in
\citet{CarrollME}] to recover the relationship of interest. For
instance, in
the epidemiological study of the effect of particulate matter pollution on
respiratory health we consider in this paper, suitable instruments could
include (i)~individual-level measurement of contaminant levels that can even
be biased and error-contaminated or (ii) incidence rates of diseases other
than the one of interest that are known to be affected by the
contaminant in
question.

Our estimation method essentially proceeds by representing each of the
unknown functions in the model by a truncated series (or a flexible
functional form) and by numerically solving for the parameter values that
best fits the observable data. Although such an approach is easy to suggest
and implement, it is a challenging task to formally establish that such a
method is guaranteed to work in general. First, there is no guarantee that
the solution (i.e., parameter values that best match the distribution
of the
observable data) is unique. Second, estimation in the presence of a number
of unknown parameters going to infinity with sample size is fraught with
convergence questions. Can the postulated series represent the solution
asymptotically? Is the parameter space too large to obtain consistency? Is
the noise associated with estimating an increasing number of parameters kept
under control?

Our solution to these problems is two-fold. First, we target the most
difficult obstacle by formally establishing identification conditions under
which the regression function and the distribution of all the unobserved
variables of the model are uniquely determined by the distribution of the
observable variables. A second important aspect of our solution to the
Berkson measurement error problem is to exploit the extensive and
well-developed literature on nonparametric sieve estimation [e.g.,
\citet
{grenandersieves}, \citet{gallantsievemle}, \citet
{Shensieve}] to formally
address the potential convergence issues that arise when nonparametric
unknowns are represented via truncated series with a number of terms that
increases with sample size. These theoretical findings are supported by a
simulation study and the usefulness of the method is illustrated with an
epidemiological application to the effect of particulate matter
pollution on
respiratory health.

\section{Model and framework}

We consider a regression model of the general form
%
%
\begin{eqnarray}
\label{eqy}
Y &=&g \bigl( X^{\ast} \bigr) +\Delta Y,
\\
\label{eqxs}
X^{\ast} &=&X+\Delta X^{\ast},
\\
\label{eqz}
Z &=&h \bigl( X^{\ast} \bigr) +\Delta Z,
\end{eqnarray}
where the function $g ( \cdot) $ is the (unknown) relationship of
interest between $Y$, the observed outcome variable and $X^{\ast}$, the
\emph{unobserved} true regressor, while $\Delta Y$ is a disturbance.
Information regarding $X^{\ast}$ is only available in the form of an
observable proxy $X$ contaminated by an error $\Delta X^{\ast}$.
Equation (%
\ref{eqz}) assumes the availability of an instrument $Z$, related to $%
X^{\ast}$ via an unknown function $h ( \cdot) $ and a
disturbance~$\Delta Z$. Our goal is to estimate the function $g ( \cdot
) $ in (\ref{eqy}) nonparametrically and without assuming that the
distribution of the measurement error $\Delta X^{\ast}$ is known. [As
by-products, we will also obtain $h ( \cdot) $ and the joint
distribution of all the unobserved variables.] To this effect, we require
the following assumptions, which are very common in the literature focusing
on nonlinear models with measurement error [e.g., \citet{CarrollME},
\citet
{wangnlmbme}, \citet{Haus1991}, \citet{fandecon2}, \citet
{li1998}, \citet
{lewbel1996}].

%
%
\begin{condition}
\label{condindep}The random variables $X$, $\Delta X^{\ast}$, $\Delta
Y$, $%
\Delta Z$ are mutually independent.
\end{condition}

Note that Assumption~\ref{condindep} implies the commonly-made
``surrogate assumption'' $f_{Y|X,X^{\ast
}} ( y|x,x^{\ast} ) =f_{Y|X^{\ast}} ( y|x^{\ast} ) $,
as can be seen by the following sequence of equalities between conditional densities: $f_{Y|X,X^{\ast
}} ( y|x,x^{\ast} ) =f_{\Delta Y|X,X^{\ast}} ( y-g (
x^{\ast} ) |x,x^{\ast} ) =f_{\Delta Y|\Delta X^{\ast},X} (
y-g ( x^{\ast} ) |x^{\ast}-x,x ) =f_{\Delta Y} (
y-g ( x^{\ast} ) ) =f_{\Delta Y|X^{\ast}} ( y-g (
x^{\ast} ) |x^{\ast} ) =f_{Y|X^{\ast}} ( y|x^{\ast} )
$.

%
%
\begin{condition}
\label{condloc}The random variables $\Delta X^{\ast}$, $\Delta Y$,
$\Delta
Z $ are centered (i.e., the model's restrictions preclude replacing
$\Delta
X^{\ast}$ by $\Delta X^{\ast}+c$ for some nonzero constant $c$, and
similarly for $\Delta Y$ and $\Delta Z$; this includes either zero mean,
zero mode or zero median, e.g.).
\end{condition}

As our approach relies on the availability of an instrument $Z$ to achieve
identification, it is instructive to provide practical examples of suitable
instruments in common settings. Although\vadjust{\goodbreak} the use of instrumental variables
has historically been more prevalent in the econometrics measurement error
literature [\citet{Haus1991,Haus1995,Neweynliv,schennachnlmeiv}],
instruments are gathering increasing interest in the statistics literature,
especially in the context of measurement error problems [see Chapter 6
entitled ``Instrumental Variables'' in \citet
{CarrollME} and the numerous references therein].

Note that instrument equation (\ref{eqz}) is entirely analogous to (\ref
{eqy}), the equation generating the main dependent variable. Hence, the
instrument is nothing but another observable
``effect'' caused by $X^{\ast}$ via a general nonlinear
relationship $h ( \cdot) $. Let us consider a few examples,
which were inspired by some of the case studies found in \citet{CarrollME},
\citet{wangnlmbme} and \citet{hyslopncme}.\vspace*{-2pt}

%
%
\begin{example}
Epidemiological studies.\vspace*{-2pt}
\end{example}

In these studies, the dependent variable $Y$ is typically a measure of the
severity of a disease or condition, while the true regressor $X^{\ast
}$ is
someone's true but unobserved exposure to some contaminant. The average
concentration $X$ of this contaminant in the region where the individual
lives is, however, observed. The error on $X$ is Berkson-type because
individual-specific $X^{\ast}$ typically randomly fluctuate around the
region average $X$. In this setup, multiple plausible instruments are
available:

\begin{longlist}[(2)]
\item[(1)] A measurement of contaminant concentration in the
individual's house
(these would be error-contaminated by classical errors, since the
concentration at a given time randomly fluctuates around the time-averaged
concentration which would be relevant for the impact on health). Thanks to
the flexibility introduced by the function $h ( \cdot) $ in (\ref
{eqz}), these measurements can even be biased. They can therefore be made
with a inexpensive method (that can be noisy and not even well-calibrated),
making it practical to use at the individual level. Hence, it is possible
to combine (i) accurate, but expensive, region averages that are not
individual-specific ($X$) and (ii) inexpensive, inaccurate
individual-specific measurements ($Z$) to obtain consistent estimates.

\item[(2)] Another plausible instrument could be a measure of the
severity of
another disease or condition that is \emph{known} to be caused by the
contaminant. The fact that it is \emph{caused by} the contaminant,
introduces an error structure which is consistent with equation (\ref{eqz}).
Other measurable effects due to the contaminant (e.g., the results of saliva
or urine tests for the presence of contaminants) could also serve as
instruments. Clearly these measurements are not units of exposure, but the
function $h ( \cdot) $ can account for this.\vspace*{-2pt}
\end{longlist}

%
%
\begin{example}
Experimental studies.\vspace*{-2pt}
\end{example}

Researchers may wish to study how an effect $Y$ (e.g., the production
of some chemical) is related to some imposed external conditions $X$
(e.g.,\vadjust{\goodbreak} oven or reactor temperature), but the true conditions $X^{\ast
}$ experienced by the sample of interest may deviate randomly from the
imposed conditions (e.g., temperature may not be completely uniform).
In this case, an instrument $Z$ could be (i) another ``effect'' (e.g.,
the amount of another chemical) that is known to be caused by $X^{\ast
}$ or (ii) a measurement of $X^{\ast}$ that is specific to the sample
of interest but that may be very noisy or even biased (e.g., it could
be an easier-to-take temperature measurement after the experiment is
completed and the sample has partly cooled down).

%
%
\begin{example}
Self-reported data.
\end{example}

\citet{hyslopncme} have argued that individuals reporting data
(e.g., their
food intake, or exercise habits) are sometimes aware of the uncertainty in
their estimates of $X^{\ast}$ and, as a result, try to report an
average $X$
over all plausible estimates consistent with the information available to
them, thus leading to Berkson-type errors, because the individuals try to
make their prediction error independent from their report. In this setting,
an instrument $Z$ could be another observable outcome variable $Z$ that is
also related to $X^{\ast}$.

\section{Identification}

We now formally state conditions under which the Berkson measurement error
model can be identified with the help of an instrument. Let $\mathcal
{Y}$, $%
\mathcal{X}$, $\mathcal{X}^{\ast}$ and $\mathcal{Z}$ denote the
supports of
the distributions of the random variables $Y$, $X$, $X^{\ast}$ and $Z$,
respectively. We consider $Y,X,X^{\ast}$ and $Z$ to be jointly continuously
distributed (with $\mathcal{Y}\subset\mathbb{R}^{n_{y}}$, $\mathcal{X}%
\subset\mathbb{R}^{n_{x}}$, $\mathcal{X}^{\ast}\subset\mathbb{R}^{n_{x}}$
and $\mathcal{Z}\subset\mathbb{R}^{n_{z}}$ with $n_{z}\geq n_{x}$).
Accordingly, we assume the following.

%
%
\begin{condition}
\label{conddens}The random variables $Y,X,X^{\ast},Z$ admit a bounded
joint density with respect to the Lebesgue measure on $\mathcal{Y}\times
\mathcal{X}\times\mathcal{X}^{\ast}\times\mathcal{Z}$. All marginal and
conditional densities are also defined and bounded.
\end{condition}

We use the notation $f_{A} ( a ) $ and $f_{A|B} ( a|b ) $ to denote the
density of the random variable $A$ and the density of $A$ conditional
on $B$, respectively. Lower case letters denote specific values of the
corresponding upper case random variables. Next, as in many treatments
of errors-in-variables models [\citet{CarrollME}, \citet
{fandecon2}, \citet{livuong1998}, \citet{li1998}, Schennach
(\citeyear{schennachnlme,schennachnlmeiv})], we require various
characteristic functions to be nonvanishing. We also place regularity
constraints on the two regression functions of the model.

%
%
\begin{condition}
\label{condinv}For all $\zeta\in\mathbb{R}^{n_{z}}$, $E [ \exp(
\mathbf{i}\zeta\cdot\Delta Z ) ] \neq0$ and for all $\xi\in
\mathbb{R}^{n_{x}}$, $E [ \exp( \mathbf{i}\xi\cdot
\Delta X^{\ast} ) ] \neq0$ (where $\mathbf{i=}\sqrt{-1}$).
\end{condition}

%
%
\begin{condition}
\label{condnodup}$g\dvtx\mathcal{X}^{\ast}\mapsto\mathcal{Y}$ and
$h\dvtx\mathcal{X%
}^{\ast}\mapsto\mathcal{Z}$ are one-to-one (but not necessarily
onto).\vadjust{\goodbreak}
\end{condition}

%
%
\begin{condition}
\label{condcont}$h$ is continuous.
\end{condition}

Assumption~\ref{condnodup} is somewhat restrictive when $X^{\ast}$ has a
dimension larger or equal to the ones of $Y$ (or $Z$). Fortunately, it is
often possible to eliminate this problem by re-defining $Y$ (and $Z$)
to be
a vector containing auxiliary variables in addition to the outcome of
interest, in order to allow for enough variation in $Y$ (and~$Z$) to
satisfy Assumption~\ref{condnodup}. Each of these additional variables need
not be part of the relationship of interest per se, but does need to be
affected by $X^{\ast}$ is some way. In that sense, such auxiliary variables
would also be a type of ``instrument.'' Our
main identification result can then be stated as follows. (Note that the
theorem also holds upon conditioning on an observed variable $W$, so that
additional, correctly measured, regressors can be straightforwardly
included.)

%
%
\begin{theorem}
\label{thid}Under Assumptions~\ref{condindep}--\ref{condcont}, given the
true observed conditional density $f_{Y,Z|X}$, the solution $ (
g,h,f_{\Delta Z},f_{\Delta Y},f_{\Delta X^{\ast}} ) $ to the
functional equation
%
%
\begin{equation}\label{eqfyxz}\quad
f_{Y,Z|X} ( y,z|x ) =\int f_{\Delta Z} \bigl( z-h \bigl(
x^{\ast
} \bigr) \bigr) f_{\Delta Y} \bigl( y-g \bigl( x^{\ast}
\bigr) \bigr) f_{\Delta X^{\ast}} \bigl( x^{\ast}-x \bigr) \,dx^{\ast}
\end{equation}
for all $y\in\mathcal{Y}$, $x\in\mathcal{X}$, $z\in\mathcal{Z}$ is
unique (up to differences on sets of null probability measure). A similar
uniqueness result holds for the solution $ ( g,h,f_{\Delta Z},\break f_{\Delta
Y}, f_{\Delta X^{\ast}},f_{X} ) $ to
%
%
\begin{eqnarray}\label{eqfyxz3}
&&
f_{Y,Z,X} ( y,z,x ) \nonumber\\[-8pt]\\[-8pt]
&&\qquad=f_{X} ( x ) \int f_{\Delta Z} \bigl(
z-h \bigl( x^{\ast} \bigr) \bigr) f_{\Delta Y} \bigl( y-g \bigl(
x^{\ast} \bigr) \bigr) f_{\Delta X^{\ast}} \bigl( x^{\ast}-x \bigr)
\,dx^{\ast}.\nonumber
\end{eqnarray}
\end{theorem}

Establishing this result demands techniques radically different from
existing treatment of Berkson error models, such as the spectral
decomposition of linear operators [see \citet{carrascoHB} for a review],
which are emerging as powerful alternatives to the ubiquitous deconvolution
techniques that are typically applied in classical measurement error
problems. The proof can be found in the \hyperref[app]{Appendix} and
can be outlined as
follows. Assumption~\ref{condindep} lets us obtain the following integral
equation relating the joint densities of the observable variables to the
joint densities of the unobservable variables:
%
%
\begin{equation}\label{eqpreindep}
f_{Y,Z|X} ( y,z|x ) =\int f_{Z|X^{\ast}} \bigl( z|x^{\ast}
\bigr) f_{Y|X^{\ast}} \bigl( y|x^{\ast} \bigr) f_{X^{\ast}|X} \bigl(
x^{\ast
}|x \bigr) \,dx^{\ast}
\end{equation}
from which equation (\ref{eqfyxz}) follows directly. Uniqueness of the
solution is then shown as follows. Equation (\ref{eqpreindep}) defines the
following operator equivalence relationship:
%
%
\begin{equation}\label{eqLeqLLL}
F_{y;Z|X}=F_{Z|X^{\ast}}D_{y;X^{\ast}}F_{X^{\ast}|X},
\end{equation}
where we have introduced the following operators:
%
%
\begin{eqnarray}
\label{eqdefop} [ F_{y;Z|X}r ] ( z ) &=& \int f_{Y,Z|X} ( y,z|x ) r (
x ) \,dx,\nonumber\\
{}[ F_{Z|X^{\ast}}r ] ( z )
&=&\int f_{Z|X^{\ast}} \bigl(
z|x^{\ast} \bigr) r \bigl( x^{\ast
} \bigr) \,dx^{\ast},
\nonumber\\
{}[ F_{Z|X}r ] ( z ) &=&\int f_{Z|X} ( z|x ) r ( x ) \,dx,\\
{}[D_{y;X^{\ast}}r ] \bigl( x^{\ast
} \bigr) &=&f_{Y|X^{\ast}} \bigl(
y|x^{\ast} \bigr) r \bigl( x^{\ast} \bigr),
\nonumber\\
{}[ F_{X^{\ast}|X}r ] \bigl( x^{\ast} \bigr) &=&\int f_{X^{\ast
}|X}
\bigl( x^{\ast}|x \bigr) r ( x ) \,dx
\nonumber
\end{eqnarray}
for some sufficiently regular but otherwise arbitrary function $r$. Note
that, in the above definitions, $y$ is viewed as a parameter (the operators
do not act on it) and that $D_{y;X^{\ast}}$ is the operator equivalent
of a
diagonal matrix. Next, we note that the equivalence
$F_{Z|X}=F_{Z|X^{\ast
}}F_{X^{\ast}|X}$ also holds [e.g., by integration of (\ref{eqLeqLLL}) over
all $y\in\mathcal{Y}$]. We can then isolate $F_{X^{\ast}|X}$
%
%
\begin{equation}\label{eqgivexsx}
F_{X^{\ast}|X}=F_{Z|X^{\ast}}^{-1}F_{Z|X}
\end{equation}
and substitute the result into (\ref{eqLeqLLL}) to yield, after
rearrangements,
%
%
\begin{equation}\label{eqdiag}
F_{y;Z|X}F_{Z|X}^{-1}=F_{Z|X^{\ast}}D_{y;X^{\ast}}F_{Z|X^{\ast}}^{-1},
\end{equation}
where all inverses can be shown to exist over suitable domains under our
assumptions. Equation (\ref{eqdiag}) states\vspace*{1pt} that the operator $%
F_{y;Z|X}F_{Z|X}^{-1}$ admits a spectral decomposition. The operator to be
``diagonalized'' is defined in terms of
observable densities, while the resulting eigenvalues $f_{Y|X^{\ast
}} ( y|x^{\ast} ) $ (contained in $D_{y;X^{\ast}}$) and
eigenfunctions $f_{Z|X^{\ast}} ( \cdot|x^{\ast} ) $ (contained
in $F_{Z|X^{\ast}}$) provide the unobserved densities of interest.

A few more steps are required to ensure uniqueness of this decomposition,
which we now briefly outline. One needs to (i) invoke a powerful uniqueness
result regarding spectral decompositions [Theorem XV 4.5 in \citet
{dunfordoper}], (ii) exploit the fact that densities integrate to one to
fix the scale of the eigenfunctions, (iii) handle degenerate
eigenvalues and
(iv) uniquely determine the ordering and indexing of the eigenvalues and
eigenfunctions. This last, and perhaps most difficult, step, addresses the
issue that both $f_{Z|X^{\ast}} ( \cdot|x^{\ast} ) $ and $%
f_{Z|X^{\ast}} ( \cdot|S ( x^{\ast} ) ) $, for some
one-to-one function $S$, are equally valid ways to state the eigenfunctions
that nevertheless result in different operators $F_{Z|X^{\ast}}$. To
resolve this ambiguity, we note that for any possible operator
$F_{Z|X^{\ast
}}$ satisfying (\ref{eqdiag}), there exist a unique corresponding
operator $%
F_{X^{\ast}|X}$, via equation (\ref{eqgivexsx}). However, only one choice
of $F_{Z|X^{\ast}}$ leads to an operator $F_{X^{\ast}|X}$ whose
kernel $%
f_{X^{\ast}|X} ( x^{\ast}|x ) $ satisfies Assumption~\ref{condloc}.
Hence, $f_{X^{\ast}|X} ( x^{\ast}|x ) $, $f_{Y|X^{\ast}} ( y|x^{\ast
} ) $ and $f_{Z|X^{\ast}} ( z|x^{\ast} ) $ are
identified, from which the functions $f_{\Delta Z}$, $f_{\Delta Y}$, $%
f_{\Delta X^{\ast}}$, $h$ and $g$ can be recovered by exploiting the
centering restrictions on $\Delta X^{\ast}$, $\Delta Y$ and $\Delta Z$.

An operator approach has recently been proposed to address certain
types of
nonclassical measurement error problems [\citet{huschennachncme}],
but under
assumptions that rule out Berkson-type measurement errors: it should be
emphasized that, despite the use of operator decomposition techniques
similar to the ones found in \citet{huschennachncme} (hereafter
HS), it is
impossible to simply use their results to identify the Berkson measurement
error model considered here, for a number of reasons. First, the key
condition (Assumption 5 in HS) that the distribution of the mismeasured
regressor $X$ given the true regressor $X^{\ast}$ is
``centered'' around $X^{\ast}$ does not hold for Berkson
errors. Consider the simple case where the Berkson measurement error is
normally distributed and so are the true and mismeasured regressors. The
distribution of $X$ given $X^{\ast}=x^{\ast}$ is a normal centered at
$%
x^{\ast}\sigma_{x}^{2}/ ( \sigma_{x}^{2}+\sigma_{\Delta x^{\ast
}}^{2} ) $. Hence, there is absolutely no reasonable measure of
location (mean, mode, median, etc.) that would yield the appropriate
centering at $x^{\ast}$ that is needed in Assumption 5 of HS. In addition,
one cannot simply replace the assumption of centering of $X$ given
$X^{\ast
} $ (as in HS) by a centering of $X^{\ast}$ given $X$ (as would be required
for Berkson errors) and hope that Theorem 1 in HS remains valid. HS exploit
the fact that, in a conditional density, there is no Jacobian term
associated with a change of variable in a conditioning variable (here $%
X^{\ast}$). However, with Berkson errors, the corresponding change of
variable would not take place in the conditional variables, and a Jacobian
term would necessarily appear, which makes the approach used in HS
fundamentally inapplicable to the Berkson case. Solving this problem
involves (i) using a different operator decomposition than in HS and (ii)
using a completely different approach for
``centering'' the mismeasured variable.\looseness=1

A referee suggested an alternative argument (formalized in the
\hyperref[app]{Appendix})
that makes a more direct connection with Theorem 1 in HS but under the
additional assumption that $Z$ and $X^{\ast}$ have the same dimension. Such
an assumption is rather restrictive because it will often result in the
assumption that $h ( \cdot) $ is one-to-one (Assumption \ref
{condnodup}) being violated. For instance, if $X^{\ast}$ is scalar and we
have access to two instruments $Z_{1}$ and $Z_{2}$ such that neither $E%
[ Z_{1}|X^{\ast} ] $ nor $E [ Z_{2}|X^{\ast} ] $ are
strictly monotone, then $h ( \cdot) $ is not one-to-one for
either instrument used in isolation. However, the mapping $X^{\ast
}\mapsto
( E [ Z_{1}|X^{\ast} ],E [ Z_{2}|X^{\ast} ]
) $ will typically be one-to-one, except for really exceptional
cases. Hence, allowing for the dimensions of $X^{\ast}$ and $Z$ to
differ is important. Nevertheless, even assuming away this problem,
such an
approach still requires a different technique for centering $X^{\ast}$ than
the one used in HS. That said, both HS and the current paper rely on
operator spectral decomposition as an alternative to conventional
convolution/deconvolution techniques, and it appears likely that these new
techniques will find applications in a number of other measurement error
models.

Observe that our identification result is also useful in a parametric and
semi-parametric context, as it provides the confidence that, under simple
conditions, the model is identified. Rank conditions that would need to be
verified on a case-by-case basis in any given parametric model are
automatically implied by our identification results in a wide class of
models. Also, although $X$ is allowed to be random throughout,
considering $%
X$ to be fixed poses no particular difficulty, since equation (\ref
{eqfyxz}%
) provides a valid conditional likelihood function in that case.

As discussed in \citet{schennachberksupp}, a number of extensions
of the
method are possible: (i) Relaxing the independence between $X$ and
$\Delta
X^{\ast}$ to allow for some heteroskedasticity in the measurement error
and (ii) combining classical and Berkson errors, a possibility considered
in, for example, \citet{carrollnevada}, \citet
{carrollberksonmix}, \citet
{stramclasberk} and \citet{hyslopncme}. It can also be shown that some
extensions are not plausible, such as assuming that both the measurement
equation (\ref{eqxs}) and the instrument equation (\ref{eqz}) have a
Berkson error structure [\citet{schennachberksupp}].

\section{Estimation}
\label{secest}

A natural way to obtain a nonparametric estimator of the model is to
substitute truncated series approximations into (\ref{eqfyxz}) or (\ref
{eqfyxz3}) for each of the unknown functions and construct a log
likelihood function to be maximized numerically with respect to all
coefficients of the series [e.g., \citet{Shensieve}]. Such sieve-based
estimators have recently found applications in a variety of measurement
error problems [e.g., \citet{Neweynliv}, \citet{mahajanmisbin},
\citet{huschennachncme}, \citet{hutwosample}, among others].
Below we
first define our estimator before establishing its consistency.

We represent the regression functions $g ( \cdot) $ and $h (
\cdot) $ as
%
%
\begin{equation}\label{expg}
\hat{m}^{ ( K_{m} ) } \bigl( x^{\ast},\beta_{m}^{ (
K_{m} ) }
\bigr) =\sum_{k=1}^{K_{m}}\beta
_{m,k}^{ (
K_{m} ) }q_{k}^{ ( K_{m} ) } \bigl(
x^{\ast} \bigr) \qquad\mbox{for }m=g,h,
\end{equation}
where $q_{k}^{ ( K_{m} ) } ( x^{\ast} ) $ is some
sequence (indexed by the truncation parameters $K_{m}$) of progressively
larger sets of basis functions indexed by $k=1,\ldots,K_{m}$ while
$\beta
_{m}^{ ( K_{m} ) }= ( \beta_{m,1}^{ ( K_{m} )
},\ldots,\beta_{m,K}^{ ( K_{m} ) } ) $ is a vector of
coefficients to be determined. The $q_{k}^{ ( K_{m} ) } (
x^{\ast} ) $ could be some power series, trigonometric series,
orthogonal polynomials, wavelets or splines, for instance. The double
indexing by $k$ and $K_{m}$ is useful to allow for splines, where changing
the number of knots modifies all the basis functions.

A similar expansion in terms of basis functions $p_{k}^{ ( K_{V} )
} ( v ) $ (with truncation parameter $K_{V}$) is used for the
density of each disturbance $V=\Delta Z,\Delta Y,\Delta X^{\ast}$,
%
%
\begin{equation}\label{expf}
\hat{f}_{V}^{ ( K_{V} ) } \bigl( v,\theta_{V}^{ ( K_{V} )}
\bigr) =\frac{1}{\theta_{V,0}^{ ( K_{v} ) }}\phi_{0} \bigl(
{v}/{\theta
_{V,0}^{ ( K_{v} ) } } \bigr) \sum_{k=1}^{K_{V}}
\theta_{V,k}^{ ( K_{V} ) }p_{k}^{ (
K_{V} ) } ( v ),
\end{equation}
where $\theta_{V}^{ ( K_{V} ) }= ( \theta_{V,0}^{ (
K_{V} ) },\ldots,\theta_{V,K}^{ ( K_{V} ) } ) $ is a
vector of coefficients to be determined, and $\phi_{0} ( \cdot) $
is a user-specified ``baseline''
function. The ``baseline'' function is
convenient to reduce the number of terms needed in the expansion, when the
approximate general shape of the density is known. It is not strictly
needed, however, and can be set to 1. Either way, the method is fully
nonparametric. A convenient choice of basis [see \citet{gallantsievemle}]
is to take $\phi_{0} ( \cdot) $ to be a Gaussian and $%
p_{k}^{ ( K_{V} ) } ( v ) =v^{k-1}$ for any $K_{V}$.

An important distinction with the functions $g ( \cdot) $ and $%
h ( \cdot) $ is that some constraints have to be imposed on the
densities. One constraint is needed to ensure centering (Assumption \ref
{condloc}),
\[
\sum_{k=1}^{K_{V}}\theta_{V,k}^{ ( K_{V} )
}C_{V,c,k}^{ ( K_{V} ) }=0,
\]
where, for some user-specified function $c_{V} ( v ) $, we define
\[
C_{V,c,k}^{ ( K_{V} ) }=\int c_{V} ( v ) \frac{1}{\theta
_{V,0}^{ ( K_{v} ) }}
\phi_{0} \biggl( \frac{v}{\theta_{V,0}^{ (
K_{v} ) }} \biggr) p_{k}^{ ( K_{V} ) }
( v ) \,dv.
\]
For instance, to impose zero mean on the disturbance $V$, let $c_{V} (
v ) =v$. To impose zero median, let $c_{V} ( v ) =\mathbf{1}%
( v\leq0 ) -1/2$, where $\mathbf{1} ( \cdot) $ denotes
an indicator function, while to impose zero mode, let $c_{V} ( v )
=-\delta^{ ( 1 ) } ( v ) $ (a delta function derivative,
in a slight abuse of notation). Another constraint is needed to ensure unit
total probability: $\sum_{k=1}^{K_{V}}\theta_{V,k}^{ ( K_{V} )
}C_{V,1,k}^{ ( K_{V} ) }=1$. Note that both types of constraints
exhibit the computationally convenient property of being linear in the
unknown coefficients.

Given the above definitions, we can define an estimator of all unknown
functions based on a sample $ ( X_{i},Y_{i},Z_{i} ) _{i=1}^{n}$ and
equation (\ref{eqfyxz}) [a corresponding estimator based on equation
(\ref
{eqfyxz3}) can be derived analogously]. Let $\hat{\beta}_{g}^{ (
K_{g} ) },\hat{\beta}_{h}^{ ( K_{g} ) },\hat{\theta}_{\Delta
X^{\ast}}^{ ( K_{V} ) },\hat{\theta}_{\Delta Y}^{ (
K_{V} ) },\hat{\theta}_{\Delta Z}^{ ( K_{V} ) }$ denote the
minimizer of the sample log likelihood
%
%
\begin{equation}\label{eqnopti}
\frac{1}{n}\sum_{i=1}^{n}\ln
\hat{f}_{Y,Z|X} ( Y_{i},Z_{i}|X_{i} ),
\end{equation}
where $\hat{f}_{Y,Z|X} ( y,z|x ) = \int\hat{f}_{\Delta Z}^{ (
K_{\Delta Z} ) } ( z-\hat{h}^{ ( K_{h} ) } ( x^{\ast
},\beta_{h}^{ ( K_{h} ) } ),\theta_{\Delta Z}^{ (
K_{\Delta Z} ) } ) \hat{f}_{\Delta Y}^{ ( K_{\Delta Y} ) }
( y-\break\hat{g}^{ (
K_{g} ) } ( x^{\ast},\beta_{g}^{ ( K_{g} ) } ),\theta_{\Delta Y}^{ (
K_{\Delta Y} ) } )$ $\hat{f}_{\Delta
X^{\ast}}^{ ( K_{\Delta X^{\ast}} ) } ( x^{\ast}-x,\theta
_{\Delta X^{\ast}}^{ ( K_{\Delta X^{\ast}} ) } ) \,dx^{\ast}$,
subject to
%
%
\begin{equation}\label{eqnlincons}
\sum_{k=1}^{K_{V}}\theta_{V,k}^{ ( K_{V} )
}C_{V,1,k}^{ ( K_{V} ) }=1
\quad\mbox{and}\quad\sum_{k=1}^{K_{V}}%
\theta
_{V,k}^{ ( K_{V} ) }C_{V,c,k}^{ ( K_{V} ) }=0
\end{equation}
for $V=\Delta Z,\Delta Y,\Delta X^{\ast}$ and subject to technical
regularity constraints to be defined below.
Estimators are then given by
%
%
\begin{eqnarray}\label{eqnest}
\hat{g} \bigl( x^{\ast} \bigr) &=&\hat{g}^{ ( K_{g} ) } \bigl(
x^{\ast};\hat{\beta}_{g}^{ ( K_{g} ) } \bigr),\qquad\hat{h} \bigl(
x^{\ast} \bigr) =\hat{h}^{ ( K_{g} ) } \bigl( x^{\ast};\hat{
\beta}%
_{h}^{ ( K_{g} ) } \bigr),
\nonumber\\[-8pt]\\[-8pt]
\hat{f}_{V} ( v ) &=&\hat{f}_{V}^{ ( K_{V} ) } \bigl(
v,%
\hat{\theta}_{V}^{ ( K_{V} ) } \bigr) \qquad\mbox{for }V=
\Delta X^{\ast
},\Delta Y,\Delta Z.
\nonumber
\end{eqnarray}
This type of estimator falls within the very general class of sieve
nonparametric maximum likelihood estimators (MLE), whose asymptotic theory
has received considerable attention over the last few decades [e.g.,
\citet
{grenandersieves}, \citet{gallantsievemle}, \citet
{Shensieve}]. Here, we
parallel the treatment of \citet{gallantsievemle} and \citet
{neweynpiv} to
establish the consistency of the above procedure. Although the consistency
of sieve-type estimators has been previously established in very general
settings under some high-level assumptions, our contribution is to provide
very primitive sufficient conditions for consistency for the class of models
considered here.

We first need to define the set in which the densities of interest reside.
The formal proof of consistency of the estimator requires this set to be
compact, although this requirement appears to have little impact in
practice. In essence, compactness is helpful to rule out very extreme but
rare events associated with very poor estimates. It is a standard regularity
condition; see, for example, \citet{gallantsievemle}, \citet
{neweynpiv},
\citet
{Neweynliv}. A well-known type of infinite-dimensional but compact sets
are those generated via boundedness and Lipschitz constraints in an $%
\mathcal{L}_{\infty}$ space. Here, we use a weighted Lipschitz constraint
in order to allow for densities supported on an unbounded set, while still
maintaining compactness (our treatment can be straightforwardly adapted to
cover the simpler case where the variables are supported on finite
intervals). Following \citet{gallantsievemle}, we enforce
restrictions that
avoid too rapid divergences in the log likelihood.

%
%
\begin{definition}
\label{defsetF}Let $\llVert f\rrVert=\sup_{v\in\mathbb{R}%
}\llvert f ( v ) \rrvert$. Let $B$ be finite and strictly
positive. Let $f_{+}^{\prime} ( v ) $ be strictly positive and
bounded function that is decreasing in~$\llvert v\rrvert$, symmetric
about $v=0$ and such that $\int_{-\infty}^{\infty}f_{+}^{\prime} (
v ) \,dv<\infty$. Let $\mathcal{S}=\{f\dvtx\mathbb{R}\mapsto[ -B,B%
] $ such that $\llvert\partial^{\lambda}f ( v ) /\partial
v^{\lambda}\rrvert\leq f_{+}^{\prime} ( v ) \}$. Let $f_{-} ( v ) $
and $f_{+} ( v ) $ be strictly positive and bounded functions with
$f_{-} ( v ) $ decreasing in $\llvert v\rrvert$ and $\int_{-\infty
}^{\infty}f_{+} ( v ) \,dv<\infty$. Let $\mathcal{F}= \{ f\in
\mathcal{S}\dvtx f_{-} ( v ) \leq f ( v ) \leq f_{+} ( v ) \} $.
\end{definition}

We also define suitable norms and sets for the regression functions. Here,
we need to allow for functions that diverge to infinity at controlled rates
toward infinite values of their argument. In analogy with any existing
global measure of expected error, we also use a norm that downweights errors
in the tails, which is consistent with the fact that the tails of a\vadjust{\goodbreak}
nonparametric regression function are always estimated with more noise,
since there are fewer datapoints there.\vspace*{-1.5pt}

%
%
\begin{definition}
\label{defsetG}Let $\omega\dvtx\mathbb{R}\mapsto\mathbb{R}^{+}$ by
some given
strictly positive, bounded and differentiable weighting function. For any
function $g\dvtx\mathbb{R}\mapsto\mathbb{R}$, let $\llVert g\rrVert
_{\omega}=\llVert\omega g\rrVert$ where $\omega g ( v )
\equiv g ( v ) \omega( v ) $. Let $\mathcal{G}= \{
g\dvtx\omega g\in\mathcal{S}$ and $\llvert g ( v ) \rrvert
\leq g_{+} ( v ) \} $ where $g_{+} ( v ) $ is a
given positive function that is increasing in $\llvert v\rrvert$ and
symmetric about $v=0$.\vspace*{-1.5pt}
\end{definition}

We can now state the regularity conditions needed.\vspace*{-1.5pt}

%
%
\begin{condition}
\label{condiid}The observed data $ ( X_{i},Y_{i},Z_{i} ) $ are
independent and identically distributed across $i=1,2,\ldots\,$.\vspace*{-1.5pt}
\end{condition}

%
%
\begin{condition}
\label{condball}We have $f_{\Delta X^{\ast}},f_{\Delta Y},f_{\Delta
Z}\in
\mathcal{F}$ and $g,h\in\mathcal{G}$.\vspace*{-1.5pt}
\end{condition}

%
%
\begin{condition}
\label{conddense}The set of functions representable as series (\ref{expf})
and (\ref{expg}) are, respectively, dense in $\mathcal{F}$ (in the norm
$%
\llVert\cdot\rrVert$) and $\mathcal{G}$ (in the norm $\llVert
\cdot\rrVert_{\omega}$).\vspace*{-1.5pt}
\end{condition}

Denseness results for numerous types of series are readily available in the
literature [e.g., \citet{neweyseries}, \citet
{gallantsievemle}]. Although
such results are sometimes phrased in a mean square-type norm rather than
the sup norm used here, Lemma~\ref{lemnorm2} below [proven in
\citet{schennachberksupp}] establishes that,
within the sets $\mathcal{F}$ and $\mathcal{G}$, denseness in a mean square
norm implies denseness in the norms we use.\vspace*{-1.5pt}

%
%
\begin{lemma}
\label{lemnorm2}Let $ \{ f_{n} \} $ be a sequence in $\mathcal{F}$.
Then $\int\llvert f_{n} ( v ) \rrvert^{2}\,dv\rightarrow0$
implies $\llVert f_{n}\rrVert\rightarrow0$ (for $\mathcal{F}$ and $%
\llVert\cdot\rrVert$ as in Definition~\ref{defsetF}).\vspace*{-1.5pt}
\end{lemma}

We also need standard boundedness and dominance conditions.\vspace*{-1.5pt}

%
%
\begin{condition}
\label{conddom}For any $x\in\mathbb{R}$, $\int( \omega( x^{\ast
} ) ) ^{-1}f_{+} ( x^{\ast}-x ) \,dx^{\ast}<\infty$ for
$\omega$ and $f_{+}$ as in Definitions~\ref{defsetG} and~\ref{defsetF},
respectively.\vspace*{-1.5pt}
\end{condition}

%
%
\begin{condition}
\label{condmomex} There exists $b>0$ such that $E [ \llvert\ln
( f_{-} ( X,Y,Z ) ) \rrvert] <\infty$, where
$f_{-} ( x,y,z ) \equiv2bf_{-} ( b ) f_{-} (
\llvert y\rrvert+ ( g_{+} ( \llvert x\rrvert
+b ) ) ) f_{-} ( \llvert z\rrvert+ (
g_{+} ( \llvert x\rrvert+b ) ) ) $ for $f_{-}$
and $g_{+}$ as in Definitions~\ref{defsetF} and~\ref{defsetG}, respectively.\vspace*{-1.5pt}
\end{condition}

We can then state our consistency result [proven in \citet
{schennachberksupp}]:\vspace*{-1.5pt}

%
%
\begin{theorem}
\label{thconsis}Under Assumptions~\ref{conddens}--\ref{condmomex}, if
$K_{V}%
\stackrel{p}{\rightarrow}\infty$, for $V=h,g,\break\Delta X^{\ast},\Delta
Y,\Delta Z$, the estimators given by (\ref{eqnest}) evaluated at the
minimizer of (\ref{eqnopti}) subject to (\ref{eqnlincons}), $\hat{f}%
_{\Delta X^{\ast}},\hat{f}_{\Delta Y},\hat{f}_{\Delta Z}\in\mathcal{F}$
and $\hat{g},\hat{h}\in\mathcal{G}$ and satisfying Assumption~\ref{conddom}
are such that $\llVert\hat{g}-g^{\ast}\rrVert_{\omega}\stackrel
{p}{%
\rightarrow}0$, $\llVert\hat{h}-h^{\ast}\rrVert_{\omega}\stackrel
{%
p}{\rightarrow}0$, $\llVert\hat{f}_{\Delta X^{\ast}}-f_{\Delta
X^{\ast
}}^{\ast}\rrVert\stackrel{p}{\rightarrow}0$, $\llVert\hat{f}%
_{\Delta Y}-f_{\Delta Y}^{\ast}\rrVert\stackrel{p}{\rightarrow}0$, $
\llVert\hat{f}_{\Delta Z}-f_{\Delta Z}^{\ast}\rrVert\stackrel{p}{%
\rightarrow}0$, where the stared quantities denote the true values [i.e.,
the unique solution to (\ref{eqfyxz})].\vadjust{\goodbreak}
\end{theorem}

The practical implementation of the above approach necessitates the
selection of the number of terms $K_{V}$ in each of the approximating
series. Theorem~\ref{thconsis} allows for a data-driven selection of
the $%
K_{V}$, since $K_{V}$ is allowed to be random. To select the $K_{V}$, one
can employ the bootstrap cross-validation model selection method based on
the Kullback--Leibler (KL) criterion, shown by \citet{laanmlecv}
to be
consistent even when the number of candidate models grows to infinity with
sample size (as it is here). In this method, a fraction $p$ of the
sample is
excluded at random and the remaining $1-p$ fraction is used to estimate the
model parameters with given numbers $ ( K_{\Delta X^{\ast}},K_{\Delta
Y},K_{\Delta Z},K_{g},K_{h} ) $ of terms in the corresponding series.
The likelihood (or KL criterion) is then evaluated using the excluded
fraction $p$ at the value of the estimated parameters found in the previous
step. The process is repeated many times with different random
partitions of
the sample into fractions $p$ and $ ( 1-p ) $, to obtain an average
KL criterion with a sufficiently small variance (which can be estimated
from the KL criterion of each random partitions). This procedure is carried
out for various trial choices of $ ( K_{\Delta X^{\ast}},K_{\Delta
Y},K_{\Delta Z},K_{g},K_{h} ) $ and the choice that yields the largest
likelihood is selected. This method is consistent asymptotically (as sample
size $n\rightarrow\infty$) as $np\rightarrow\infty$ and
$p\rightarrow0$
and under some mild technical regularity conditions stated in
\citet
{laanmlecv}.

Our nonparametric approach nests parametric and semiparametric models. These
subcases can be easily implemented by replacing some, or all, of the
nonparametric series approximations by suitable parametric models. It is
possible to obtain convergence rates and limiting distribution results,
along the lines of \citet{Shensieve} or \citet
{huschennachncme}, although
we do not do so here due to space limitations [stating suitable regularity
conditions, even in high-level form, is rather involved, as seen in the
supplementary material of \citet{huschennachncme}, which covers a related
but different measurement error model]. It is, however, important to point
out one important property. Sieve nonparametric MLE is optimal in the
following sense: under suitable regularity conditions, any sufficiently
regular semiparametric functional of the nonparametric sieve MLE estimates
is asymptotically normal and root $n$ consistent and reaches the
semiparametric efficiency bound for that functional; see Theorem 4 in
\citet
{Shensieve}. This notion of optimality is a natural nonparametric
generalization of the well-known efficiency of parametric maximum likelihood.

\section{Simulations study}
\label{secsimul}

We now investigate the practical performance and feasibility of the
proposed estimator via a simulation example purposely chosen to be a
difficult case. The data is generated as follows. The distribution of
$X$ is a uniform distribution over $ [ -1,1 ] $ (implying a standard
deviation of $0.58$). We consider a thick-tailed $t$ distribution with
6 degrees of freedom scaled by $0.5$ as the distribution of $\Delta
X^{\ast}$. The standard deviation of the error $\Delta X^{\ast}$ is
almost identical to the one of the ``signal'' $X$, thus making this
estimation problem exceedingly difficult. The distribution
of $%
\Delta Y$ is a logistic scaled by $0.125$ while the distribution of
$\Delta
Z $ is a $t$ distribution with $6$ degrees of freedom scaled by $0.25$. The
regression function has the form\looseness=1
%
%
\begin{equation}
g \bigl( x^{\ast} \bigr) =\bigl\llvert x^{\ast}\bigr\rrvert
x^{\ast},
\end{equation}\looseness=0
which is only finitely many times differentiable, thus limiting the
convergence rate of its series estimator in the measurement-error-robust
estimator (the naive estimator would be less affected since it would
``see'' a smoothed version of this
function). The instrument equation has a specification that is strictly
convex and therefore tends to exacerbate the bias in many nonparametric
estimators,
\[
h \bigl( x^{\ast} \bigr) =\ln\bigl( 1+\exp\bigl( 2x^{\ast}
\bigr) \bigr).
\]

A total of $100$ independent samples, each containing $500$ observations,
were generated as above and fed into our estimator. For estimation purposes,
the functions $g ( \cdot) $ and $h ( \cdot) $ are both
represented by polynomials while the densities of $\Delta X^{\ast}$, $%
\Delta Y$ and $\Delta Z$ are represented by a Gaussian multiplied by a
polynomial [following \citet{gallantsievemle}, who establish that these
choices satisfy a suitable denseness condition]. The Gaussian is centered
at the origin, but its width is left as a parameter to be estimated. Note
that the functional forms considered are not trivially nested within the
space spanned by the truncated sieve approximation. This was an intentional
choice aimed at properly accounting for the nonparametric nature of the
problem (in which the researcher never has the fortune of selecting a
truncated sieve fitting the true model exactly).

The integral in equation (\ref{eqfyxz}) is evaluated numerically by
discretizing the integral as a sum over the range $ [ -3,3 ] $ in
intervals of $0.05$. Naive least-squares estimators ignoring measurement
error (i.e., least-squares regressions of $Y$ on $X$ and of $Z$ on $X$) were
used as a starting point for the numerical sieve optimization of the
$g$ and
$h$ functions, while the variances of the corresponding residuals were used
to construct an initial Gaussian guess for the optimization of all the error
distributions. The simplex method due to \citet{nelderamoeba} (also
known as
``amoeba'') was used to carry out the
numerical optimization of the log likelihood (\ref{eqnopti}) with
respect to
all the parameters $\theta_{V}^{ ( K_{V} ) }$ for $V=\Delta
X^{\ast},\Delta Y,\Delta Z$ and $\beta_{m}^{ ( K_{m} ) }$ for $%
m=g,h$ simultaneously. The constraints that the estimated densities and
regression functions lie, respectively, in the sets $\mathcal{F}$ and $%
\mathcal{G}$ of the form given in Definitions~\ref{defsetF} and~\ref{defsetG}
are implied by bounds on the magnitude of the sieve coefficients $\theta
_{V,k}^{ ( K_{V} ) }$ and $\beta_{m,k}^{ ( K_{m} ) }$ in
(\ref{expf}) and (\ref{expg}). Such constraints are easy to impose within
the simplex optimization method: parameter changes that would yield
violations of the bounds are simply rejected (effectively assigned an
``infinite'' value)---the simplex
optimization method easily accommodates such extreme behavior in the
objective function, since it does not rely on derivatives. However, we
found that these constraints are rarely binding in practice, unless the
number of terms $K_{V}$ in the expansions is large [\citet{gallantsievemle}
reports a similar observation]. Such large values of $K_{V}$ tend to be
naturally ruled out via our data-driven selection method of the number of
terms.

To select the number of terms in the approximating series for a given
sample, we use the ``bootstrap
cross-validation'' method described in Section~\ref{secest}
with a fraction $p=1/8$ and $100$ bootstrap replications. Trial values of
the number of free parameters (not counting parameters uniquely determined
by zero mean and unit area constraints) in the series representing $%
f_{\Delta X^{\ast}},f_{\Delta Y},f_{\Delta Z}$ each span the set $ \{
1,2,3,4 \} $ while for $g$, $h$ each span the set $ \{
4,5,6,7 \} $. The optimal numbers of parameters (kept constant during
the replications) were found to be $f_{\Delta X^{\ast}}:3$; $f_{\Delta
Y}:3 $; $f_{\Delta Z}:3$; $g:6$; $h:6$.

\begin{figure}

\includegraphics{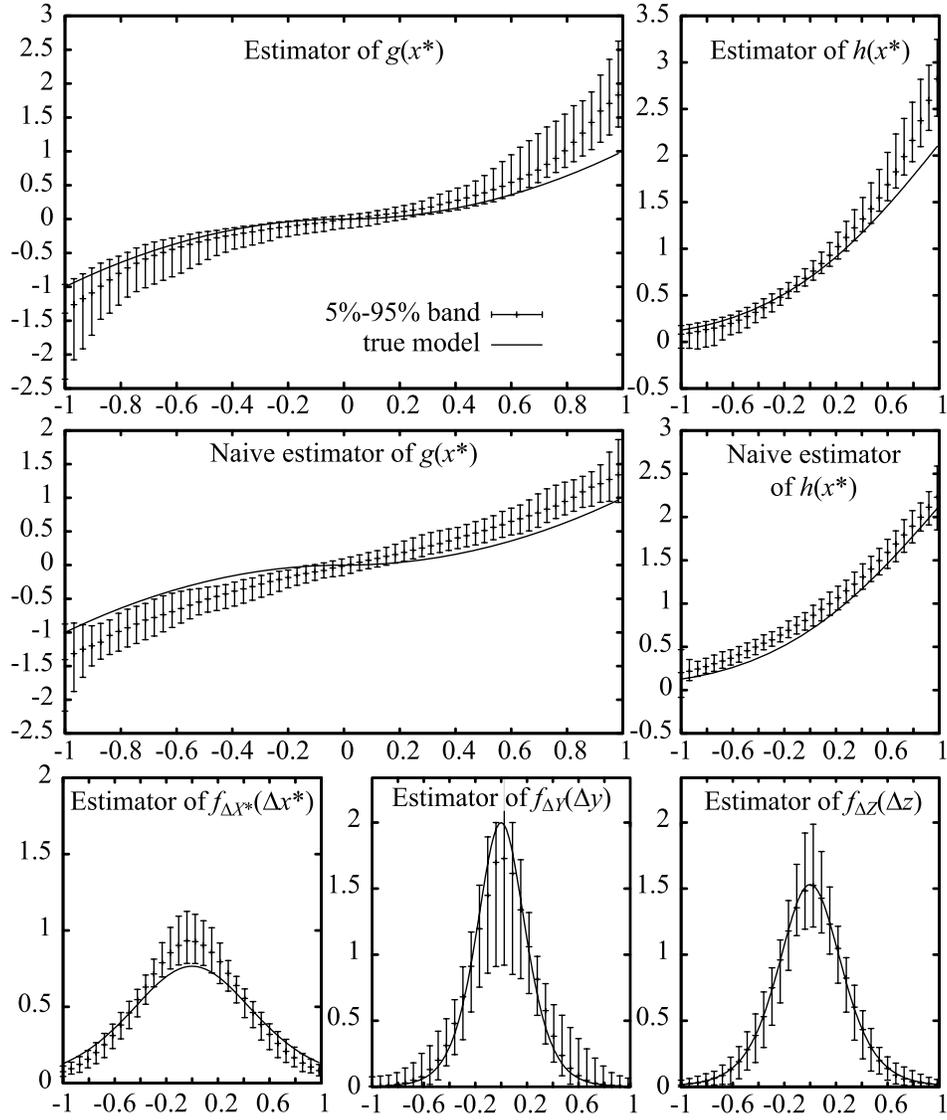}

\caption{Simulation study of the practical performance of the proposed
measurement-error-robust
estimator in comparison
with a ``naive'' nonparametric polynomial series least-square estimator
that ignores
the presence of measurement error.
In each plot, the pointwise 90\% confidence band of the estimator
simulated over
100 replications is shown as error bars.}\label{figberksim}\vspace*{-3pt}
\end{figure}

Figure~\ref{figberksim} summarizes the result of these simulations,
where a
naive nonparametric series least-squares estimator ignoring measurement
error (i.e., least-square regressions of $Y$ on $X$ and of $Z$ on $X$) with
the same number of sieve terms is also shown for comparison. The reliability
of the method can be appreciated by noting how closely the median of the
replicated measurement-error-robust estimates matches the true model, while
the naive estimator ignoring the presence of measurement error is
considerably more biased, even missing the fact that the true regression
function is nearly flat in the middle section and instead producing a very
misleading linear shape despite the strong nonlinearity of the true model.
In fact, unlike the proposed estimator, the naive estimator is so
significantly biased that any type of hypothesis test based on it would
exhibit completely misleading confidence levels: the true model curves
(for $%
g$ and $h$) almost always lies beyond the 95\% or 5\% percentiles of the
estimator distribution.

%
%

Overall, the proposed measurement-error-robust estimator exhibits low
variability and low bias at the reasonable sample size of $500$. The bias
is not exactly zero in a finite sample because our estimator is a nonlinear
functional of sample averages and because the sieve approximation
necessarily has a limited accuracy in a finite sample. Nevertheless, the
fact that our estimator performs so well in the presence of measurement
error of such large magnitude is a strong indication of its practical
usefulness. This behavior is not specific to this model---we have tested
the method in other simulation settings; see \citet{schennachberksupp}.

\section{Application}
\label{secapp}

Numerous studies have sought to quantify the effect of air pollution on
respiratory health [e.g., \citet{dockerypollmort}]. Specifically, there
is a growing concern regarding the effect of small particulate matter
[\citet{popepartair}, \citet{sametpartair}]. A key difficulty
with such
studies is that air quality monitors are not necessarily located near
the subjects being affected by air pollution, implying that the main
regressor of interest is mismeasured.

Our approach to this question relies on very comprehensive country-wide data
collected by Environment Protection Agency (EPA) and the Center for Disease
Control (CDC) in the United States. Pollution levels are taken from EPAs
Monitor Values Report---Criteria Air Pollutants database for year 2005.
EPAs data provides point measurements of the particulate matter levels (we
focus on so-called 95th percentile level of PM2.5 particles, those having
less than 2.5 micrometers in diameter) at various monitoring stations
throughout the United States, from which we construct state-averaged
pollution levels (our $X$ variable, measured in $\mu$g of particles
per m$%
^{3}$). We do so because pollution data is only available for a small
fraction of counties, and even where it is available, the nature of its
measurement error is complex (it could be a mixture of classical and Berkson
errors). By constructing state-level averages, we average out the
randomness in monitor measurements while leaving the randomness in the
individual exposure untouched, thus obtaining a valid Berkson
error-contaminated estimate of the pollution level experienced by
individuals from each state, whether they live in a county with a monitoring
station or not. Each individual faces an exposure equal to the state average
plus an unknown random noise due to his/her precise geographic whereabouts
and lifestyle.

Health data is obtained from the publicly available ``CDC Wonder''
data\-base entitled ``Mortality---underlying cause of death'' for year
2005. To measure respiratory health, we use data on causes of death,
which offers the advantage that it is very comprehensive and accurate
(medical professionals are required to collect it and there is no
reliance on voluntary surveys). One limit to the completeness of the
data is that, for some counties, the data is ``suppressed'' (for
privacy reasons) or labeled as ``unreliable'' by the CDC and were
therefore omitted from our sample. Our dependent variable of interest
($Y$) is the rate (per 10,000) of death due to ``chronic lower
respiratory diseases'' (e.g., asthma, bronchitis, emphysema), while our
instrument ($Z$) is the rate (per 10,000) of death resulting from
``lung diseases due to external agents'' (e.g., pneumoconiosis due to
organic or inorganic dust, coalworker's pneumoconiosis). The rationale
is to use, as an instrument, a~variable that is clearly expected to be
affected by pollution levels. This variable indirectly provides
information regarding the true level of pollution, so that the effect
of pollution (if any) on the variable of interest can be more
accurately assessed. We employ county-level data on causes of death
because they are readily available without concerns for patient privacy
issues. Moreover, the CDC provides age-corrected death rates, thus
correcting for demographic differences between counties. We construct
our sample by matching mortality data via counties and matching
pollution data via states, resulting in 1305 observations over as many
counties and covering all 51 states. A~limitation of our approach is
that it does not control for other possible confounding effects, for
example, if the proportion of smokers differs between industrial and
nonindustrial cities. However, such a limitation is common in studies
of this kind [as noted in \citet{dockerypollmort}].

We use the same types of sieves and computational methods as in the
simulation example and select the number of terms using the
``bootstrap cross-validation'' method
described in Section~\ref{secest} with a fraction $p=1/8$ and $100$
bootstrap replications. Trial values of the number of free parameters in
the series representing $f_{\Delta X^{\ast}},f_{\Delta Y},f_{\Delta Z}$
span the range $ \{ 1,2,3 \} $ while trial values of the number of
terms in the series representing $g$ and $h$ span the range $ \{
2,3,4 \} $ (increasing any one of the $K_{V}$ beyond that range
resulted in clearly worse performances). The optimal numbers of free
parameters (not counting parameters uniquely determined by zero mean and
unit area constraints) were found to be $f_{\Delta X^{\ast}}:2$; $%
f_{\Delta Y}:3$; $f_{\Delta Z}:1$; $g:4$; $h:3$. Pointwise 90\% confidence
bands around the nonparametric estimates were obtained using the standard
bootstrap [see, e.g., \citet{ginezinnboot} for general conditions
justifying its use] with 100 replications.

%
%
\begin{figure}[b]\vspace*{6pt}

\includegraphics{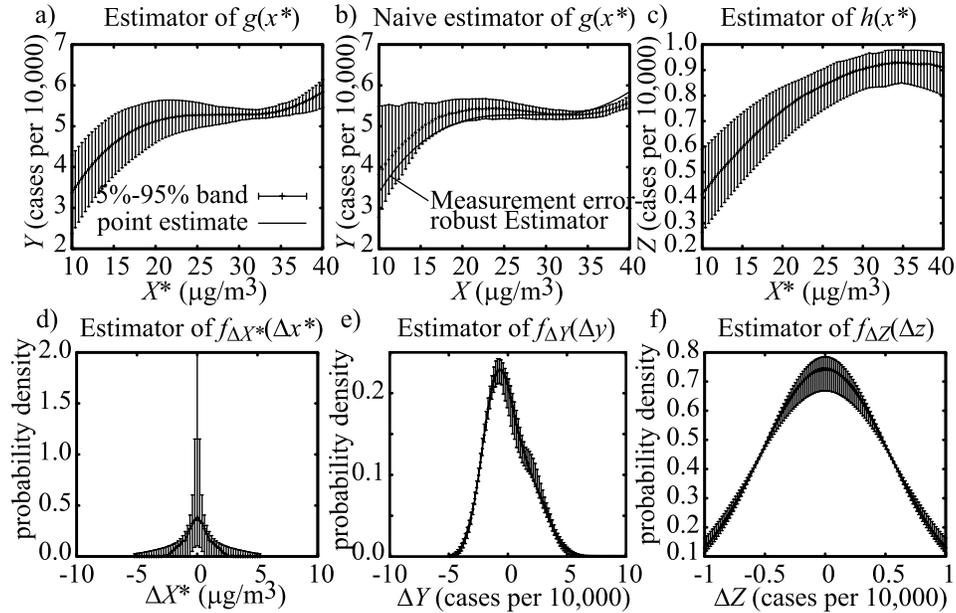}

\caption{Application of the proposed estimator to an epidemiologic
example (see text for a description of the variables and the estimated
functions). In each plot, the estimator is shown as a solid line while
the error bars indicate the pointwise 90\% confidence bands. In \textup{(b)},
the ``naive'' estimator is a nonparametric polynomial series
least-square estimator that ignores the presence of measurement error.
The estimator in \textup{(a)} is shown on the plot \textup{(b)} for comparison.}
\label{figapp}
\end{figure}

Results are shown in Figure~\ref{figapp}. A few observations are in order.
First, our measurement error-robust estimator is perfectly able to
detect a
clear monotone relationship between $Y$ and $X^{\ast}$ and between $Z$
and $%
X^{\ast}$ with useful confidence bands, despite the use of a fully
nonparametric approach. Second, although the distribution of the measurement
error is difficult to estimate (as reflected by the wide confidence bands),
the impact of this uncertainty on the main function of interest [$g (
x^{\ast} ) $] is fortunately very limited. The 90\% confidence bands
indicate that the presence of substantial measurement error is consistent
with the data: the measurement error is of the order of 10 $\mu$g/m$^{3}$,
whereas the observed $X$ roughly ranges from 10 to 40 $\mu$g/m$^{3}$.
Third, the distribution of $\Delta Y$ exhibits nonnegligible asymmetry, thus
illustrating the drawbacks of methods merely assuming normality of all the
error terms. In contrast, the distributions of $\Delta X^{\ast}$ and $%
\Delta Z$ are apparently very close to symmetric (this is a conclusion of
the formal model selection procedure, not an assumption).

For comparison purposes, we also naively regress the dependent
variables ($Y$~or~$Z$) on the mismeasured regressor $X$ using a conventional least
squares (thereby neglecting measurement error) with a polynomial
specification with the same number of terms as our Berkson model. A first
troubling observation from this exercise [see Figure~\ref{figapp}(b)] is
that the naive estimate of $g ( x^{\ast} ) $ is not monotone,
although in the region where it is unexpectedly decreasing, the confidence
bands do not rule out a constant response. Second, it is perhaps
counter-intuitive that the confidence bands for the naive estimator are
sometimes larger than the corresponding bands for the measurement
error-robust estimator. This is a consequence of the fact that correcting
for Berkson errors amounts to an operation akin to convolution (rather than
deconvolution, as in classical measurement errors). Unlike deconvolution,
convolution is a noise-reducing operation, effectively averaging
observations of $Y$ over a wide range of values of $X$ to yield an estimate
the expected value of $Y$ given a specific value of $X$. This
phenomenon is
probably also responsible for the more reasonable (i.e., increasing) behavior
of the response for the measurement error-robust estimate. Finally, the
measurement error-robust regression function often lies at or beyond the
95\% or 5\% percentiles of the naive estimator distribution; see Figure
\ref{figapp}(b). This implies that the level of any statistical test
would be
severely biased. For instance, the confidence bands of the naive estimator
would reject our best estimate of $g ( x^{\ast} ) $ obtained with
the measurement-error robust procedure.

In summary, this application example serves to illustrate that ignoring
Berkson errors can be seriously misleading in nonlinear settings. Not only
is the shape of the estimated response considerably affected, but
statistical inferences based on a measurement error-blind method would be
seriously biased. This application example also shows that our fully
nonparametric and measurement error-robust method works well at sample sizes
typically available in real data sets, without assuming the knowledge
of the
distribution of the measurement error.\vspace*{-2pt}

%
\begin{appendix}\label{app}
\section*{Appendix: Proofs}\vspace*{-2pt}

Let $\mathcal{L}_{1}^{b} ( \mathcal{D} ) $ with $\mathcal{D}\subset
$ $\mathbb{R}^{n_{0}}$ for some $n_{0}$ denote the set of all bounded
functions in $\mathcal{L}_{1} ( \mathcal{D} ) $ endowed with the\vadjust{\goodbreak}
usual $\mathcal{L}_{1}$ norm. Also, whenever we state an equality between
functions in $\mathcal{L}_{1}^{b} ( \mathcal{D} ) $, we mean that
their difference is zero in the $\mathcal{L}_{1}$ norm.

We provide two proofs of Theorem 1. The first one, suggested by a referee,
relies on the additional assumptions that (i) $Z$ and $X^{\ast}$ have the
same dimension and (ii) $h$ and its inverse are differentiable. Assumption
(i) makes Assumption~\ref{condnodup} unlikely to hold, but enables a
somewhat direct application of Theorem 1 in \citet
{huschennachncme}. The
second proof relaxes those assumptions. It borrows some of the operator
techniques from \citet{huschennachncme}, yet requires considerable changes
in the approach---we focus here on the aspects of the proof that differ.

\begin{pf*}{Proof Theorem \protect\ref{thid} (simple special case)}
Let variables from \citet{huschennachncme} be denoted by the corresponding
uppercase letter with tildes and make the following assignments: $ (
\tilde{X}^{\ast},\tilde{X},\tilde{Y},\tilde{Z} ) = ( h (
X^{\ast} ),Z,Y,X ) $. We now verify the 5 assumptions of Theorem
1 in \citet{huschennachncme}.

To verify Assumption 1, we observe that the densities of $ ( \tilde{X}%
^{\ast},\tilde{X},\tilde{Y},\tilde{Z} ) $ and $ ( X^{\ast
},Z,Y,X ) $ are related through:%
$f_{\tilde{X}^{\ast},\tilde{X},\tilde{Y},\tilde{Z}} ( \tilde{x}^{\ast
},%
\tilde{x},\tilde{y},\tilde{z} ) =f_{X^{\ast},Z,Y,X} ( h^{-1} (
\tilde{x}^{\ast} ),\break\tilde{x},\tilde{y},\tilde{z} ) \llvert
{\partial h^{-1} ( \tilde{x}^{\ast} ) }/{\partial\tilde{x}%
^{\ast\prime}}\rrvert$
where the density $f_{X^{\ast},Z,Y,X}$ exists by Assumption~\ref{conddens},
and $h^{-1} ( \tilde{x}^{\ast} ) $ exists by Assumption \ref
{condnodup}. The Jacobian $\partial h^{-1} ( \tilde{x}^{\ast} )
/\partial\tilde{x}^{\ast\prime}$ matrix is only defined if $X^{\ast}$
and $Z$ (and therefore $\tilde{X}^{\ast}$) have the same dimension and is
finite and nonsingular under the assumption that $h$ and its inverse are
differentiable. A~similar argument can be used for marginals and
conditional distributions.

To verify Assumption 2, we note that our model can be written in terms of
tilded variables as
%
%
\begin{eqnarray}
\label{eqtilde1}
\tilde{Y} &=&Y=g \bigl( h^{-1} \bigl( \tilde{X}^{\ast} \bigr)
\bigr) +\Delta Y,
\\
\label{eqtilde2}
\tilde{Z} &=&X=h^{-1} \bigl( \tilde{X}^{\ast} \bigr) -\Delta
X^{\ast},
\\
\label{eqtilde3}
\tilde{X} &=&Z=\tilde{X}^{\ast}+\Delta Z.
\end{eqnarray}
To verify Assumption 2(i), we write
\begin{eqnarray*}
f_{\tilde{Y}|\tilde{X},\tilde{X}^{\ast},\tilde{Z}} \bigl(
\tilde{y}|\tilde{x},\tilde{x}^{\ast},
\tilde{z} \bigr) &=&f_{Y|Z,X^{\ast},X} \bigl( \tilde{y}|%
\tilde{x},h^{-1} \bigl( \tilde{x}^{\ast} \bigr),\tilde{z} \bigr)
\\
&=&f_{\Delta Y|\Delta Z,\Delta X^{\ast},X} \bigl( \tilde{y}-g \bigl(
h^{-1} \bigl(
\tilde{x}^{\ast} \bigr) \bigr) |\tilde{x}-\tilde{x}^{\ast
},h^{-1}
\bigl( \tilde{x}^{\ast} \bigr) -\tilde{z},\tilde{z} \bigr)
\\
&=&f_{\Delta Y} \bigl( \tilde{y}-g \bigl( h^{-1} \bigl(
\tilde{x}^{\ast
} \bigr) \bigr) \bigr)
\\
&=&f_{Y|\tilde{X}^{\ast}} \bigl( \tilde{y}|\tilde{x}^{\ast} \bigr)
=f_{%
\tilde{Y}|\tilde{X}^{\ast}} \bigl( \tilde{y}|\tilde{x}^{\ast} \bigr),
\end{eqnarray*}
where we have used, in turn, (i) the equality $ ( \tilde{X}^{\ast},%
\tilde{X},\tilde{Y},\tilde{Z} ) = ( h ( X^{\ast} ),Z,Y,X ) $ and the
fact that changes of variables in the conditioning variables do not
introduce Jacobian terms, (ii) the fact that conditioning on $Z,X^{\ast
},X$ is equivalent to conditioning on $\Delta Z,\Delta X^{\ast
},X$, (iii) Assumption~\ref{condindep}, (iv) the relationship between $%
\Delta Y$ and $Y$ via (\ref{eqtilde1}) and (v) the equality $Y=\tilde{Y}$.

To verify Assumption 2(ii), we similarly write
\begin{eqnarray*}
f_{\tilde{X}|\tilde{X}^{\ast},\tilde{Z}} \bigl( \tilde{x}|\tilde
{x}^{\ast},%
\tilde{z}
\bigr) &=&f_{Z|X^{\ast},X} \bigl( \tilde{x}|h^{-1} \bigl(
\tilde{x}%
^{\ast} \bigr),\tilde{z} \bigr)
\\
&=&f_{\Delta Z|\Delta X^{\ast},X} \bigl( \tilde{x}-\tilde{x}^{\ast
}|h^{-1}
\bigl( \tilde{x}^{\ast} \bigr) -\tilde{z},\tilde{z} \bigr)
\\
&=&f_{\Delta Z} \bigl( \tilde{x}-\tilde{x}^{\ast} \bigr) =
f_{Z|\tilde{X}^{\ast}} \bigl( \tilde{x}|\tilde{x}^{\ast} \bigr)
=f_{%
\tilde{X}|\tilde{X}^{\ast}} \bigl( \tilde{x}|\tilde{x}^{\ast} \bigr).
\end{eqnarray*}

Assumption 3 is implied by Assumptions~\ref{conddens},~\ref{condindep},
\ref{condinv},~\ref{condnodup},~\ref{condcont} and Lem\-ma~\ref{leminj} below.

Assumption 4 requires that $f_{\tilde{Y}|\tilde{X}^{\ast}} ( \tilde{y}|
\tilde{x}_{1}^{\ast} ) \not=f_{\tilde{Y}|\tilde{X}^{\ast}} (
\tilde{y}|\tilde{x}_{2}^{\ast} ) $ for $\tilde{x}_{1}^{\ast}\not=%
\tilde{x}_{2}^{\ast}$. This can be verified as follows:
\begin{eqnarray*}
f_{\tilde{Y}|\tilde{X}^{\ast}} \bigl( \tilde{y}|\tilde{x}_{1}^{\ast}
\bigr)
&=&f_{\Delta Y|\tilde{X}^{\ast}} \bigl( \tilde{y}-g \bigl( h^{-1} \bigl(
\tilde{x}_{1}^{\ast} \bigr) \bigr) |\tilde{x}_{1}^{\ast}
\bigr)
\\
&=&f_{\Delta Y} \bigl( \tilde{y}-g \bigl( h^{-1} \bigl(
\tilde{x}_{1}^{\ast
} \bigr) \bigr) \bigr)
\\
&\neq&f_{\Delta Y} \bigl( \tilde{y}-g \bigl( h^{-1} \bigl(
\tilde{x}_{2}^{\ast
} \bigr) \bigr) \bigr) =f_{\tilde{Y}|\tilde{X}^{\ast}}
\bigl( \tilde{y}|%
\tilde{x}_{2}^{\ast} \bigr)
\end{eqnarray*}
by invoking (i) the definition of $\Delta Y$, (ii) independence of
$\Delta Y$
from $X^{\ast}$ (and therefore $\tilde{X}^{\ast}$), (iii) the fact
that $%
\tilde{x}_{1}^{\ast}\neq\tilde{x}_{2}^{\ast}$ implies $g (
h^{-1} ( \tilde{x}_{1}^{\ast} ) ) \not=g ( h^{-1} (
\tilde{x}_{2}^{\ast} ) ) $ since $g ( \cdot) $ and $%
h ( \cdot) $ are one-to-one by Assumption~\ref{condnodup} and so
is $g ( h^{-1} ( \cdot) ) $.

Assumption 5 is trivially satisfied, by equation (\ref{eqtilde3}).

Theorem 1 in \citet{huschennachncme} then allows us to conclude
that the
joint distribution of $ ( h ( X^{\ast} ),X,Y,Z ) $ is
identified. However, in order to identify the distribution of $ (
X^{\ast},X,Y,Z ) $, we need to identify $h ( \cdot) $. To
this effect, we note that, conditional on $X=x$, the fluctuations in $%
\tilde{X}^{\ast}$ are entirely caused by fluctuations in $\Delta
X^{\ast}$
by equation (\ref{eqtilde2}). Moreover, $\Delta X^{\ast}$ is
independent from $X$, hence
%
%
\begin{equation}\label{eqfJ}
f_{\tilde{X}^{\ast}|X} \bigl( \tilde{x}^{\ast}|x \bigr) =f_{\Delta
X^{\ast
}}
\bigl( h^{-1} \bigl( \tilde{x}^{\ast} \bigr) -x \bigr) \biggl
\llvert\frac{%
\partial h^{-1} ( \tilde{x}^{\ast} ) }{\partial\tilde{x}^{\ast
\prime}}\biggr\rrvert,
\end{equation}
where the left-hand side was previously identified and where the Jacobian
term is well defined by Assumptions~\ref{condnodup} and the assumed
differentiability of $h^{-1} ( \tilde{x}^{\ast} ) $. The Jacobian
can be identified by integrating (\ref{eqfJ}) with respect to $x^{\ast
}$ to
yield $\int f_{\tilde{X}^{\ast}|X} ( \tilde{x}^{\ast}|x ) \,dx
=\llvert\frac{\partial h^{-1} (
\tilde{x}^{\ast} ) }{\partial\tilde{x}^{\ast\prime}}\rrvert$.
By varying $x$ while keeping $\tilde{x}^{\ast}$ fixed in equation (\ref
{eqfJ}), we can identify the density $f_{\Delta X^{\ast}}$ up to a
shift of
$h^{-1} ( \tilde{x}^{\ast} ) $. Assumption~\ref{condloc}, pins
down what the shift should be, so that $h^{-1} ( \tilde{x}^{\ast
} ) $ is identified for any given~$\tilde{x}^{\ast}$. Since $h (
\cdot) $ is one-to-one by Assumption~\ref{condnodup}, $h^{-1} (
\cdot) $ uniquely determines $h ( \cdot) $. Hence, the
joint distribution of $ ( X^{\ast},X,Y,Z ) $ is identified.
Finally, noting that $f_{Y|X^{\ast}} ( y|x^{\ast} ) =f_{\Delta
Y} ( y-g ( x^{\ast} ) ) $ (by Assumption~\ref{condindep}), then
establishes the identification of $g ( x^{\ast} ) $ with
the help of Assumption~\ref{condloc}.
\end{pf*}

\begin{pf*}{Proof of Theorem \protect\ref{thid} (general case)}
This proof borrows some of the operator techniques from \citet
{huschennachncme}, and we focus here on the aspects of the proof that differ.

The definition of marginal and conditional densities in combination with
Assumption~\ref{condindep} lead to the following sequence of equalities:
\begin{eqnarray*}
&&
f_{Y,Z|X} ( y,z|x ) \\
&&\qquad=\int f_{Y|X^{\ast},Z,X} \bigl( y|x^{\ast
},z,x
\bigr) f_{X^{\ast},Z|X} \bigl( x^{\ast},z|x \bigr) \,dx^{\ast}
\\
&&\qquad=\int f_{\Delta Y|X^{\ast},\Delta Z,\Delta X^{\ast}} \bigl( y-g
\bigl( x^{\ast} \bigr)
|x^{\ast},z-h \bigl( x^{\ast} \bigr),x^{\ast}-x \bigr)
f_{X^{\ast},Z|X} \bigl( x^{\ast},z|x \bigr) \,dx^{\ast}
\\
&&\qquad=\int f_{\Delta Y} \bigl( y-g \bigl( x^{\ast} \bigr) \bigr)
f_{X^{\ast
},Z|X} \bigl( x^{\ast},z|x \bigr) \,dx^{\ast}
\\
&&\qquad=\int f_{\Delta Y} \bigl( y-g \bigl( x^{\ast} \bigr) \bigr)
f_{Z|X^{\ast
},X} \bigl( z|x^{\ast},x \bigr) f_{X^{\ast}|X} \bigl(
x^{\ast}|x \bigr) \,dx^{\ast}
\\
&&\qquad=\int f_{\Delta Y} \bigl( y-g \bigl( x^{\ast} \bigr) \bigr)
f_{\Delta
Z|X^{\ast},\Delta X^{\ast}} \bigl( z-h \bigl( x^{\ast} \bigr)
|x^{\ast
},x^{\ast}-x
\bigr)\\
&&\qquad\quad\hspace*{8.5pt}{}\times f_{\Delta X^{\ast}|X} \bigl( x^{\ast}-x|x \bigr) \,dx^{\ast}
\\
&&\qquad=\int f_{\Delta Y} \bigl( y-g \bigl( x^{\ast} \bigr) \bigr)
f_{\Delta
Z} \bigl( z-h \bigl( x^{\ast} \bigr) \bigr)
f_{\Delta X^{\ast}} \bigl( x^{\ast}-x \bigr) \,dx^{\ast}
\end{eqnarray*}
or, equivalently,
%
%
\begin{equation}\label{eqpreop1}
f_{Y,Z|X} ( y,z|x ) =\int f_{Z|X^{\ast}} \bigl( z|x^{\ast}
\bigr) f_{Y|X^{\ast}} \bigl( y|x^{\ast} \bigr) f_{X^{\ast}|X} \bigl(
x^{\ast
}|x \bigr) \,dx^{\ast}.
\end{equation}
As in \citet{huschennachncme}, this integral equation can be
written more
conveniently as an operator equivalence relation
%
%
\begin{equation}\label{eqop1}
F_{y;Z|X}=F_{Z|X^{\ast}}D_{y;X^{\ast}}F_{X^{\ast}|X}
\end{equation}
by introducing the operators defined in equation (\ref{eqdefop}), which are
acting on an arbitrary $r\in\mathcal{%
L}_{1}^{b} ( \mathcal{X} ) $ [or $r\in\mathcal{L}_{1}^{b} (
\mathcal{X}^{\ast} ) $].%

Similarly, one can show that
%
%
\begin{equation}\label{eqfzxsidx}
f_{Z|X} ( z|x ) =\int f_{Z|X^{\ast}} \bigl( z|x^{\ast} \bigr)
f_{X^{\ast}|X} \bigl( x^{\ast}|x \bigr) \,dx^{\ast}
\end{equation}
and thus $F_{Z|X}=F_{Z|X^{\ast}}F_{X^{\ast}|X}$.
By Assumptions~\ref{conddens},~\ref{condindep},~\ref{condinv}, \ref
{condnodup},~\ref{condcont} and Lem\-ma~\ref{leminj} below, we know that $
F_{Z|X^{\ast}}$ admits an inverse on the range of $F_{Z|X^{\ast}}$ (and
therefore the range of $F_{Z|X}$), and we can write
%
%
\begin{equation}\label{eqop2i}
F_{X^{\ast}|X}=F_{Z|X^{\ast}}^{-1}F_{Z|X}.
\end{equation}
Substituting (\ref{eqop2i}) into (\ref{eqop1}), we obtain
\[
F_{y;Z|X}=F_{Z|X^{\ast}}D_{y;X^{\ast}}F_{Z|X^{\ast}}^{-1}F_{Z|X}.
\]
By Assumptions~\ref{conddens},~\ref{condindep},~\ref{condinv}, \ref
{condnodup},~\ref{condcont} and Lemma~\ref{leminj} below again, $F_{Z|X}$
admits an inverse. Moreover, by Lemma 1 in \citet
{huschennachncme}, the
domain of $F_{Z|X}^{-1}$ is dense in $\mathcal{L}_{1}^{b} ( \mathcal{Z}%
) $, and we can then write
%
%
\begin{equation}\label{eqdiagpr}
F_{y;Z|X}F_{Z|X}^{-1}=F_{Z|X^{\ast}}D_{y;X^{\ast}}F_{Z|X^{\ast}}^{-1}.
\end{equation}
Equation (\ref{eqdiagpr}) states that the operator $F_{y;Z|X}F_{Z|X}^{-1}$
admits a spectral decomposition, where the eigenvalues are given by the
$%
f_{Y|X^{\ast}} ( y|x^{\ast} ) $ for $x^{\ast}\in\mathcal{X}%
^{\ast}$ (for a fixed $y$) defining the operator $D_{y;X^{\ast}}$ while
the eigenfunctions are the functions $f_{Z|X^{\ast}} ( \cdot|x^{\ast
} ) $ for $x^{\ast}\in\mathcal{X}^{\ast}$ defining the kernel of the
operator $F_{Z|X^{\ast}}$. As usual, the knowledge of a linear operator
[e.g., $F_{Z|X}$] only determines the value of its kernel [e.g., $%
f_{Z|X} ( z|x ) $] everywhere except on a set of null Lebesgue
measure. The resulting equivalence class exactly matches the usual
equivalence class for probability densities with respect to the Lebesgue
measure, so identifiability of the model is not affected.

The operator to be diagonalized is entirely defined in terms of observable
densities while the decomposition provides the unobserved densities of
interest. To ensure uniqueness of this decomposition, we employ four
techniques. First, a powerful result from spectral analysis [Theorem XV 4.5
in \citet{dunfordoper}] ensures uniqueness up to some normalizations.
Second, the a priori arbitrary scale of the eigenfunctions is fixed
by the requirement that densities must integrate to one. Third, to avoid
any ambiguity in the definition of the eigenfunctions when degenerate
eigenvalues are present, we use Assumption~\ref{condnodup} and the fact that
the eigenfunctions [which do not depend on $y$, unlike the eigenvalues $
f_{y|x^{\ast}} ( y|x^{\ast} ) $] must be consistent across
different values of the dependent variable~$y$. These three steps are
described in detail in \citet{huschennachncme} and are not
repeated here.

The fourth step [which differs from the approach taken in \citet
{huschennachncme}] is to rule out that the eigenvalues $f_{y;X^{\ast
}} ( y,x^{\ast} ) $ and eigenfunctions $f_{Z|X^{\ast}} (
\cdot|x^{\ast} ) $ could be indexed by a different variable without
affecting the operator $F_{y;Z|X}F_{Z|X}^{-1}$. (This issue is
analogous to
the nonunique ordering of the eigenvalues and eigenvectors in matrix
diagonalization.) Suppose that the eigenfunctions can be indexed by another
value, that is, they are given by $f_{Z|\tilde{X}^{\ast}} ( \cdot
|\tilde{%
x}^{\ast} ) $ where $\tilde{x}^{\ast}$ is another variable related to
$x^{\ast}$ through $x^{\ast}=S ( \tilde{x}^{\ast} ) $ for some
one-to-one function~$S$.\setcounter{footnote}{1}\footnote{%
Note that $S ( \cdot) $ is also measurable, for otherwise
$X^{\ast} \equiv S( \tilde{X}^{\ast} )$ would not be a proper random
variable.}
Under this alternative indexing, all the assumptions of the
original model must still hold with $x^{\ast}$ replaced by $\tilde
{x}^{\ast
}$, so a relationship similar to (\ref{eqfzxsidx}) would still have to hold,
for the same observed $f_{Z|X} ( z|x ) $
%
%
\begin{equation}\label{eqaltfzxs}
f_{Z|X} ( z|x ) =\int f_{Z|\tilde{X}^{\ast}} \bigl( z|\tilde{x}%
^{\ast} \bigr) f_{\tilde{X}^{\ast}|X} \bigl( \tilde{x}^{\ast}|x \bigr)
\,d%
\tilde{x}^{\ast}
\end{equation}
or, in operator notation, $F_{Z|X}=F_{Z|\tilde{X}^{\ast}}F_{\tilde
{X}^{\ast}|X}$.

In order for $f_{Z|\tilde{X}^{\ast}} ( z|\tilde{x}^{\ast} ) $ to
be a valid alternative density, it must satisfy the same assumptions (and
their implications) as $f_{Z|X^{\ast}} ( z|x^{\ast} ) $. In
particular, the fact that $F_{Z|X^{\ast}}$ is invertible (established
above via Lemma~\ref{leminj}) must also hold for $F_{Z|\tilde{X}^{\ast}}$.
Hence, for any alternative $F_{Z|\tilde{X}^{\ast}}$, there is a unique
corresponding $F_{\tilde{X}^{\ast}|X}$, given by $F_{\tilde{X}^{\ast
}|X}=F_{Z|\tilde{X}^{\ast}}^{-1}F_{Z|X}$. We can find a more explicit
expression for $f_{\tilde{X}^{\ast}|X} ( \tilde{x}^{\ast}|x ) $
as follows. First note that we trivially have that $f_{Z|\tilde{X}^{\ast
}} ( z|\tilde{x}^{\ast} ) =f_{Z|X^{\ast}} ( z|S ( \tilde{x%
}^{\ast} ) ) $ since $x^{\ast}=S ( \tilde{x}^{\ast} )
$ and $S$ is one-to-one. By performing the change of variable $x^{\ast
}=S ( \tilde{x}^{\ast} ) $ in (\ref{eqfzxsidx}), we obtain
\[
f_{Z|X} ( z|x ) =\int f_{Z|X^{\ast}} \bigl( z|S \bigl(
\tilde{x}%
^{\ast} \bigr) \bigr) f_{X^{\ast}|X} \bigl( S \bigl(
\tilde{x}^{\ast
} \bigr) |x \bigr) \,d\mu\bigl( \tilde{x}^{\ast}
\bigr),
\]
where the measure $\mu$ is defined, via $\mu( \mathcal{A} )
=\lambda( S^{-1} ( \mathcal{A} ) ) $ for any measurable
set~$\mathcal{A}$, where $\lambda$ denotes the Lebesgue measure and $%
S^{-1} ( \mathcal{A} ) \equiv\{ \tilde{x}^{\ast}\in\mathcal{%
A}:S ( \tilde{x}^{\ast} ) =x^{\ast} \} $. From this
we can conclude the equality between the two following measures:
%
%
\begin{equation}\label{eqdS1}
f_{\tilde{X}^{\ast}|X} \bigl( \tilde{x}^{\ast}|x \bigr) \,d\tilde
{x}^{\ast
}=f_{X^{\ast}|X}
\bigl( S \bigl( \tilde{x}^{\ast} \bigr) |x \bigr) \,d\mu\bigl(
\tilde{x}^{\ast} \bigr)
\end{equation}
by comparison to equation (\ref{eqaltfzxs}) and the uniqueness of the
measure\break $f_{\tilde{X}^{\ast}|X} ( \tilde{x}^{\ast}| x ) \,d\tilde{x}%
^{\ast}$ due to the injectivity of the $F_{Z|\tilde{X}^{\ast}}$ operator,
shown in Lem\-ma~\ref{leminj} in the general case where the domain of
$F_{Z|%
\tilde{X}^{\ast}}$ could include finite signed measures. We will now show
that $f_{\tilde{X}^{\ast}|X} ( \tilde{x}^{\ast}|x ) $ necessarily
violates Assumption~\ref{condloc} (with $\Delta X^{\ast}$ replaced by $
\Delta\tilde{X}^{\ast}\equiv\tilde{X}^{\ast}-X$), unless $S ( \cdot
) $ is the identity function.\footnote{%
Some of the steps below were inspired by comments from an anonymous referee.}

Since $\Delta X^{\ast}=X^{\ast}-X$ with $\Delta X^{\ast}$ independent
from $X$, we have $f_{X^{\ast}|X} ( x^{\ast}|x ) =f_{\Delta
X^{\ast}} ( x^{\ast}-x ) $ and by a similar reasoning $f_{\tilde{%
X}^{\ast}|X} ( \tilde{x}^{\ast}|x ) =f_{\Delta\tilde{X}^{\ast
}} ( \tilde{x}^{\ast}-x ) $ with $\Delta\tilde{X}^{\ast}\equiv
\tilde{X}^{\ast}-X$. Equation (\ref{eqdS1}) then becomes
%
%
\begin{equation}\label{eqmueq}
f_{\Delta\tilde{X}^{\ast}} \bigl( \tilde{x}^{\ast}-x \bigr) \,d\tilde
{x}%
^{\ast}=f_{\Delta X^{\ast}} \bigl( S \bigl( \tilde{x}^{\ast} \bigr) -x
\bigr) \,d\mu\bigl( \tilde{x}^{\ast} \bigr).
\end{equation}
Now, for a given $x$, consider Radom--Nikodym derivative of $f_{\Delta
\tilde{%
X}^{\ast}} ( \tilde{x}^{\ast}-x ) \,d\tilde{x}^{\ast}$ with
respect to the Lebesgue measure $d\tilde{x}^{\ast}$, which is, by
definition (almost everywhere) equal to $f_{\Delta\tilde{X}^{\ast
}} ( \tilde{x}^{\ast}-x ) $, a bounded function by Assumption \ref
{conddens}. By equation (\ref{eqmueq}), the existence of the Radom--Nikodym
derivative of the left-hand side implies the existence of the same
Radom--Nikodym derivative on the right-hand side, and we can write
%
%
\begin{equation}\label{eqallf}
f_{\Delta\tilde{X}^{\ast}} \bigl( \tilde{x}^{\ast}-x \bigr)
=f_{\Delta
X^{\ast}}
\bigl( S \bigl( \tilde{x}^{\ast} \bigr) -x \bigr) \,\frac{d\mu
( \tilde{x}^{\ast} ) }{d\tilde{x}^{\ast}}
\end{equation}
almost everywhere. Integrating both sides of the equation over all $x\in
\mathcal{X}$, we obtain (after noting that points where the equality may
fail have null measure and therefore do not contribute to the integral),
$1=1\,\frac{d\mu( \tilde{x}^{\ast} ) }{d\tilde{x}^{\ast}}$,
since densities integrate to $1$, which implies that $d\mu( \tilde{x}%
^{\ast} ) /d\tilde{x}^{\ast}=1$, that is, $\mu$ is also the Lebesgue
measure. It follows from (\ref{eqallf}) that, almost everywhere
\[
f_{\Delta\tilde{X}^{\ast}} \bigl( \tilde{x}^{\ast}-x \bigr)
=f_{\Delta
X^{\ast}}
\bigl( S \bigl( \tilde{x}^{\ast} \bigr) -x \bigr).
\]
In order for Assumption~\ref{condloc} to hold for both $\Delta\tilde{X}
^{\ast}$ and $\Delta X^{\ast}$, we must have that $f_{\Delta\tilde{X}
^{\ast}} ( \tilde{x}^{\ast}-x ) $, when viewed as a function of $%
\tilde{x}^{\ast}$ for any given $x$, is centered at $\tilde{x}^{\ast}=x$,
and we must simultaneously have that $f_{\Delta X^{\ast}} ( x^{\ast
}-x ) =f_{\Delta X^{\ast}} ( S ( \tilde{x}^{\ast} )
-x ) $, when viewed as a function of $x^{\ast}$ for any given $x$, is
centered at $x^{\ast}=x$, that is, $S ( \tilde{x}^{\ast} ) =x$. The
two statements are only compatible if $\tilde{x}^{\ast}=S ( \tilde{x}%
^{\ast} ) $. Thus, there cannot exist two distinct but
observationally equivalent parametrization of the eigenvalues/eigenfunctions.

Hence we have shown, through equation (\ref{eqdiagpr}), that the unobserved
functions $f_{Y|X^{\ast}} ( y|x^{\ast} ) $ and $f_{Z|X^{\ast
}} ( \cdot|x^{\ast} ) $ are uniquely determined (up to an
equivalence class of functions differing at most on a set of null Lebesgue
measure) by the observed function $f_{Y,Z|X} ( y,z|x ) $. Next,
equation (\ref{eqop2i}) implies that $f_{X^{\ast}|X} ( x^{\ast
}|x ) $ is uniquely determined as well.

Once $f_{Y|X^{\ast}} ( y|x^{\ast} ) $ and $f_{Z|X^{\ast}} (
z|x^{\ast} ) $ are known, the functions $g ( x^{\ast} ) $
and $h ( x^{\ast} ) $ can be identified by exploiting the
centering restrictions on $\Delta Y$, $\Delta X^{\ast}$ and $\Delta Z$,
for example, $g ( x^{\ast} ) =\int yf_{Y|X^{\ast}} ( y|x^{\ast
} ) \,dy$ if $\Delta Y$ is assumed to have zero mean. Next, $f_{\Delta
Y} ( \Delta y ) $ can be straightforwardly identified, for example, $%
f_{\Delta Y} ( \Delta y ) =f_{Y|X^{\ast}} ( g ( x^{\ast
} ) +\Delta y|x^{\ast} ) $ for any $x^{\ast}\in\mathcal{X}%
^{\ast}$. Similar arguments yield $h ( x^{\ast} ) $ and $%
f_{\Delta Z} ( \Delta z ) $ from $f_{Z|X^{\ast}} ( z|x^{\ast
} ) $ as well as $f_{\Delta X^{\ast}} ( \Delta x^{\ast} ) $
from $f_{X^{\ast}|X} ( x^{\ast}|x ) $. It follows that equation (%
\ref{eqfyxz}) has a unique solution. The second conclusion of the theorem
then follows from the fact that both $f_{Y,Z|X} ( y,z|x ) $ and $%
f_{X} ( x ) $ are uniquely determined (except perhaps on a set of
null Lebesgue measure) from $f_{Y,Z,X} ( y,z,x ) $.
\end{pf*}

The following lemma is closely related to Proposition 2.4 in \citet
{dhaultcomp}. It is different in terms of the spaces the operators can act
on and more general in terms of the possible dimensionalities of the random
variables involved.

%
%
\begin{lemma}
\label{leminj}Let $X,X^{\ast}$ and $Z$ be generated by equations (\ref
{eqxs}%
) and (\ref{eqz}). Let $\mathcal{S} ( \mathcal{T} ) $ be the set of
finite signed measures on a given set $\mathcal{T=X},\mathcal{X}^{\ast}$
or~$\mathcal{Z}$ [and note that\vadjust{\goodbreak} $\mathcal{S} ( \mathcal{T} ) $
includes $\mathcal{L}_{1}^{b} ( \mathcal{T} ) $ as a special case,
in the sense that for any function in $r\in\mathcal{L}_{1}^{b} (
\mathcal{T} ) $, there is a corresponding measure $R\in\mathcal{S}%
( \mathcal{T} ) $ whose Radom--Nikodym derivative with respect to
the Lebesgue measure is $r$]. Under Assumptions~\ref{condindep}, \ref
{conddens},~\ref{condinv},~\ref{condnodup} and~\ref{condcont}, the
operators $F_{X^{\ast}|X}\dvtx\mathcal{S} ( \mathcal{X} ) \mapsto
\mathcal{L}_{1}^{b} ( \mathcal{X}^{\ast} ) $, $F_{Z|X^{\ast}}\dvtx%
\mathcal{S} ( \mathcal{X}^{\ast} ) \mapsto\mathcal{L}%
_{1}^{b} ( \mathcal{Z} ) $ and $F_{Z|X}\dvtx\mathcal{S} ( \mathcal{X%
} ) \mapsto\mathcal{L}_{1}^{b} ( \mathcal{Z} ) $, defined in (%
\ref{eqdefop}), are injective mappings.
\end{lemma}

\begin{pf}
First, one can verify that $R\in\mathcal{S} ( \mathcal{X} ) $
implies that $F_{X^{\ast}|X}R\in\mathcal{L}_{1}^{b} ( \mathcal{X}%
^{\ast} ) $ and similarly for $F_{Z|X^{\ast}}$ and $F_{Z|X}$, since
the (conditional) densities involving variables $X^{\ast},X$ and $Z$ are
bounded by Assumption~\ref{conddens} and are absolutely integrable. We now
verify injectivity of $F_{Z|X^{\ast}}$.

By Assumptions~\ref{condindep},~\ref{conddens} and equation (\ref
{eqz}), we
have, for any $R\in\mathcal{S} ( \mathcal{X}^{\ast} ) $,
\[
[ F_{Z|X^{\ast}}R ] ( z ) = \int f_{Z|X^{\ast}} \bigl( z|x^{\ast}
\bigr) \,dR \bigl( x^{\ast} \bigr)
= \int f_{\Delta Z} \bigl( z-h \bigl( x^{\ast} \bigr) \bigr) \,dR
\bigl( x^{\ast} \bigr).
\]
Next, let $\tilde{R}$ denote the signed measure assigning, to any measurable
set $\mathcal{A}\subseteq\mathbb{R}^{n_{z}}$, the value%
$\tilde{R} ( \mathcal{A} ) =\int1 ( h ( x^{\ast} )
\in\mathcal{A} ) \,dR ( x^{\ast} )$
and note that $\tilde{R}$ is a finite signed measure since $R ( x^{\ast
} ) $ is. Then, we can express $F_{Z|X^{\ast}}R$ as
%
%
\begin{equation}\label{eqconvol1}
[ F_{Z|X^{\ast}}R ] ( z ) =\int f_{\Delta Z} \bigl( z-%
\tilde{x}^{\ast} \bigr) \,d\tilde{R} \bigl( \tilde{x}^{\ast} \bigr),
\end{equation}
that is,\vspace*{1pt} a convolution between the probability measure of $\Delta Z$
(represented by its Lebesgue density) and the signed measure
$\tilde{R}$; see Chapter 5 in \citet{raoasympexp}. By the convolution
theorem for signed measures [Theorem 5.1(iii) in \citet{raoasympexp}],
one can convert
the convolution (\ref{eqconvol1}) into a product of Fourier transforms,%
\footnote{%
Note that the Fourier transforms involved are all continuous functions
because the original functions (or measures) are absolutely integrable (or
finite), hence ``almost everywhere''
qualifications do not apply to them.}
\[
\sigma( \zeta) =\phi_{\Delta Z} ( \zeta) \rho( \zeta)
\]
where $\sigma( \zeta) \equiv\int[ F_{Z|X^{\ast}}R ]
( z ) e^{\mathbf{i}\zeta z}\,dz$, $\phi_{\Delta Z} ( \zeta
) \equiv E [ e^{\mathbf{i}\zeta Z} ] $ and $\rho( \zeta
) \equiv\int e^{\mathbf{i}\zeta z}\,d\tilde{R} ( z ) $. Since $%
\phi_{\Delta Z} ( \zeta) $, the characteristic function of $%
\Delta Z$, is nonvanishing by Assumption~\ref{condinv}, we can isolate
$\rho
( \zeta) $ as
\[
\rho( \zeta) =\sigma( \zeta) /\phi_{\Delta
Z} ( \zeta).
\]
Since there is a one-to-one mapping between finite signed measures and their
Fourier transforms [by Theorem 5.1(i) in \citet{raoasympexp}],
$\tilde{R}$
can be recovered as the unique signed measure whose Fourier transform
is $%
\rho( \zeta) $. We now show that the signed measure $\tilde{R}$
uniquely determines the measure $R$.

Let $\mathcal{A}_{\mathcal{B}}=\bigcup_{x^{\ast}\in\mathcal{B}} \{
h ( x^{\ast} ) \} $ for any measurable $\mathcal{B}\subseteq
\mathbb{R}^{n_{x}}$, and note that $\mathcal{A}_{\mathcal{B}}$ is also
measurable since $h$ is continuous by\vadjust{\goodbreak} Assumption~\ref{condcont}. Then
observe that by Assumption~\ref{condnodup}, $h ( x^{\ast} ) \in
\mathcal{A}_{\mathcal{B}}$ if and only if $x^{\ast}\in\mathcal{B}$,
and we have
\[
\tilde{R} ( \mathcal{A}_{\mathcal{B}} ) = \int1 \bigl( h \bigl(
x^{\ast} \bigr) \in\mathcal{A}_{\mathcal{B}} \bigr) \,dR \bigl(
x^{\ast
} \bigr) = \int1 \bigl( x^{\ast}\in\mathcal{B} \bigr) \,dR
\bigl( x^{\ast} \bigr).
\]
Since $\mathcal{B}$ is arbitrary, the knowledge of $\tilde{R} ( \mathcal
{%
A}_{\mathcal{B}} ) $ uniquely determines the value assigned to any
measurable set by the signed measure $R$.

Injectivity of $F_{X^{\ast}|X}$ is a special case of the above derivation
(with $Z,X^{\ast}$ replaced by $X^{\ast},X$), in which $h$ is the
identity function. Finally, injectivity of $F_{Z|X}$ is implied by the
injectivity of $%
F_{Z|X^{\ast}}$ and $F_{X^{\ast}|X}$, since $F_{Z|X}=F_{Z|X^{\ast
}}F_{X^{\ast}|X}$ by Assumption~\ref{condindep} and equations (\ref
{eqxs}) and (\ref{eqz}).
\end{pf}
\end{appendix}


\begin{supplement}
\stitle{Supplementary material to ``Regressions with Berkson errors in
covari\-ates---A nonparametric approach''}
\slink[doi]{10.1214/13-AOS1122SUPP} 
\sdatatype{.pdf}
\sfilename{aos1122\_supp.pdf}
\sdescription{The supplementary material provides (i) a proof of
consistency of the proposed estimator, (ii) additional simulation
results and (iii) various extensions of the method, including the
weakening of some of full independence assumptions to conditional
independence and handling the simultaneous presence of classical and
Berkson errors.}
\end{supplement}

%

\printaddresses


\begin{thebibliography}{39}

\bibitem[\protect\citeauthoryear{Berkson}{1950}]{berskonme}
%
\begin{barticle}[author]
\bauthor{\bsnm{Berkson},~\bfnm{J.}\binits{J.}}
(\byear{1950}).
\btitle{Are there two regressions?}
\bjournal{J. Amer. Statist. Assoc.}
\bvolume{45}
\bpages{164--180}.
\bptok{imsref}%
\end{barticle}
%
\endbibitem

\bibitem[\protect\citeauthoryear{Bhattacharya and Rao}{2010}]{raoasympexp}
%
\begin{bbook}[author]
\bauthor{\bsnm{Bhattacharya},~\bfnm{R.~N.}\binits{R.~N.}} \AND
\bauthor{\bsnm{Rao},~\bfnm{R.~R.}\binits{R.~R.}}
(\byear{2010}).
\btitle{Normal Approximation and Asymptotic Expansions}.
\bpublisher{SIAM}, \blocation{Philadelphia}.
\bptok{imsref}%
\end{bbook}
%
\endbibitem

\bibitem[\protect\citeauthoryear{Carrasco, Florens and
Renault}{2005}]{carrascoHB}
%
\begin{bincollection}[author]
\bauthor{\bsnm{Carrasco},~\bfnm{M.}\binits{M.}},
\bauthor{\bsnm{Florens},~\bfnm{J.~P.}\binits{J.~P.}} \AND
\bauthor{\bsnm{Renault},~\bfnm{E.}\binits{E.}}
(\byear{2005}).
\btitle{Linear inverse problems and structural econometrics: Estimation based
on spectral decomposition and regularization}.
In \bbooktitle{Handbook of Econometrics, Vol. 6}.
\bpublisher{Elsevier}, \blocation{Amsterdam}.
\bptok{imsref}%
\end{bincollection}
%
\endbibitem

\bibitem[\protect\citeauthoryear{Carroll, Chen and Hu}{2010}]{hutwosample}
%
\begin{barticle}[mr]
\bauthor{\bsnm{Carroll},~\bfnm{Raymond~J.}\binits{R.~J.}},
\bauthor{\bsnm{Chen},~\bfnm{Xiaohong}\binits{X.}} \AND
\bauthor{\bsnm{Hu},~\bfnm{Yingyao}\binits{Y.}}
(\byear{2010}).
\btitle{Identification and estimation of nonlinear models using two samples
with nonclassical measurement errors}.
\bjournal{J. Nonparametr. Stat.}
\bvolume{22}
\bpages{379--399}.
\bid{doi={10.1080/10485250902874688}, issn={1048-5252}, mr={2662599}}
\bptok{imsref}%
\end{barticle}
%
\endbibitem

\bibitem[\protect\citeauthoryear{Carroll, Delaigle and
Hall}{2007}]{carrollberksonmix}
%
\begin{barticle}[mr]
\bauthor{\bsnm{Carroll},~\bfnm{Raymond~J.}\binits{R.~J.}},
\bauthor{\bsnm{Delaigle},~\bfnm{Aurore}\binits{A.}} \AND
\bauthor{\bsnm{Hall},~\bfnm{Peter}\binits{P.}}
(\byear{2007}).
\btitle{Non-parametric regression estimation from data contaminated by a
mixture of {B}erkson and classical errors}.
\bjournal{J. R. Stat. Soc. Ser. B Stat. Methodol.}
\bvolume{69}
\bpages{859--878}.
\bid{doi={10.1111/j.1467-9868.2007.00614.x}, issn={1369-7412}, mr={2368574}}
\bptok{imsref}%
\end{barticle}
%
\endbibitem

\bibitem[\protect\citeauthoryear{Carroll et~al.}{2006}]{CarrollME}
%
\begin{bbook}[mr]
\bauthor{\bsnm{Carroll},~\bfnm{Raymond~J.}\binits{R.~J.}},
\bauthor{\bsnm{Ruppert},~\bfnm{David}\binits{D.}},
\bauthor{\bsnm{Stefanski},~\bfnm{Leonard~A.}\binits{L.~A.}} \AND
\bauthor{\bsnm{Crainiceanu},~\bfnm{Ciprian~M.}\binits{C.~M.}}
(\byear{2006}).
\btitle{Measurement Error in Nonlinear Models}.
\bpublisher{Chapman \& Hall/CRC}, \blocation{Boca Raton, FL}.
\bptok{imsref}%
\end{bbook}
%
\endbibitem

\bibitem[\protect\citeauthoryear{Delaigle, Hall and
Qiu}{2006}]{delaiglenpberkson}
%
\begin{barticle}[mr]
\bauthor{\bsnm{Delaigle},~\bfnm{Aurore}\binits{A.}},
\bauthor{\bsnm{Hall},~\bfnm{Peter}\binits{P.}} \AND
\bauthor{\bsnm{Qiu},~\bfnm{Peihua}\binits{P.}}
(\byear{2006}).
\btitle{Nonparametric methods for solving the {B}erkson errors-in-variables
problem}.
\bjournal{J. R. Stat. Soc. Ser. B Stat. Methodol.}
\bvolume{68}
\bpages{201--220}.
\bid{doi={10.1111/j.1467-9868.2006.00540.x}, issn={1369-7412}, mr={2188982}}
\bptok{imsref}%
\end{barticle}
%
\endbibitem

\bibitem[\protect\citeauthoryear{D'Haultfoeuille}{2011}]{dhaultcomp}
%
\begin{barticle}[mr]
\bauthor{\bsnm{D'Haultfoeuille},~\bfnm{Xavier}\binits{X.}}
(\byear{2011}).
\btitle{On the completeness condition in nonparametric instrumental problems}.
\bjournal{Econometric Theory}
\bvolume{27}
\bpages{460--471}.
\bid{doi={10.1017/S0266466610000368}, issn={0266-4666}, mr={2806256}}
\bptok{imsref}%
\end{barticle}
%
\endbibitem

\bibitem[\protect\citeauthoryear{Dockery et~al.}{1993}]{dockerypollmort}
%
\begin{barticle}[author]
\bauthor{\bsnm{Dockery},~\bfnm{D.~W.}\binits{D.~W.}},
\bauthor{\bsnm{Pope},~\bfnm{C.~A.}\binits{C.~A.}},
\bauthor{\bsnm{Xu},~\bfnm{X.~P.}\binits{X.~P.}} \betal{et~al.}
(\byear{1993}).
\btitle{An association between air-pollution and mortality in 6 United-States
cities}.
\bjournal{New England Journal of Medicine}
\bvolume{329}
\bpages{1753--1759}.
\bptok{imsref}%
\end{barticle}
%
\endbibitem

\bibitem[\protect\citeauthoryear{Dunford and Schwartz}{1971}]{dunfordoper}
%
\begin{bbook}[author]
\bauthor{\bsnm{Dunford},~\bfnm{N.}\binits{N.}} \AND
\bauthor{\bsnm{Schwartz},~\bfnm{J.~T.}\binits{J.~T.}}
(\byear{1971}).
\btitle{Linear Operators}.
\bpublisher{Wiley}, \blocation{New York}.
\bptok{imsref}%
\end{bbook}
%
\endbibitem

\bibitem[\protect\citeauthoryear{Fan and Truong}{1993}]{fandecon2}
%
\begin{barticle}[mr]
\bauthor{\bsnm{Fan},~\bfnm{Jianqing}\binits{J.}} \AND
\bauthor{\bsnm{Truong},~\bfnm{Young~K.}\binits{Y.~K.}}
(\byear{1993}).
\btitle{Nonparametric regression with errors in variables}.
\bjournal{Ann. Statist.}
\bvolume{21}
\bpages{1900--1925}.
\bid{doi={10.1214/aos/1176349402}, issn={0090-5364}, mr={1245773}}
\bptok{imsref}%
\end{barticle}
%
\endbibitem

\bibitem[\protect\citeauthoryear{Gallant and Nychka}{1987}]{gallantsievemle}
%
\begin{barticle}[mr]
\bauthor{\bsnm{Gallant},~\bfnm{A.~Ronald}\binits{A.~R.}} \AND
\bauthor{\bsnm{Nychka},~\bfnm{Douglas~W.}\binits{D.~W.}}
(\byear{1987}).
\btitle{Semi-nonparametric maximum likelihood estimation}.
\bjournal{Econometrica}
\bvolume{55}
\bpages{363--390}.
\bid{doi={10.2307/1913241}, issn={0012-9682}, mr={0882100}}
\bptok{imsref}%
\end{barticle}
%
\endbibitem

\bibitem[\protect\citeauthoryear{Gin{\'e} and Zinn}{1990}]{ginezinnboot}
%
\begin{barticle}[mr]
\bauthor{\bsnm{Gin{\'e}},~\bfnm{Evarist}\binits{E.}} \AND
\bauthor{\bsnm{Zinn},~\bfnm{Joel}\binits{J.}}
(\byear{1990}).
\btitle{Bootstrapping general empirical measures}.
\bjournal{Ann. Probab.}
\bvolume{18}
\bpages{851--869}.
\bid{issn={0091-1798}, mr={1055437}}
\bptok{imsref}%
\end{barticle}
%
\endbibitem

\bibitem[\protect\citeauthoryear{Grenander}{1981}]{grenandersieves}
%
\begin{bbook}[mr]
\bauthor{\bsnm{Grenander},~\bfnm{Ulf}\binits{U.}}
(\byear{1981}).
\btitle{Abstract Inference}.
\bpublisher{Wiley}, \blocation{New York}.
\bid{mr={0599175}}
\bptok{imsref}%
\end{bbook}
%
\endbibitem

\bibitem[\protect\citeauthoryear{Hausman, Newey and Powell}{1995}]{Haus1995}
%
\begin{barticle}[mr]
\bauthor{\bsnm{Hausman},~\bfnm{J.~A.}\binits{J.~A.}},
\bauthor{\bsnm{Newey},~\bfnm{W.~K.}\binits{W.~K.}} \AND
\bauthor{\bsnm{Powell},~\bfnm{J.~L.}\binits{J.~L.}}
(\byear{1995}).
\btitle{Nonlinear errors in variables: Estimation of some {E}ngel curves}.
\bjournal{J. Econometrics}
\bvolume{65}
\bpages{205--233}.
\bid{doi={10.1016/0304-4076(94)01602-V}, issn={0304-4076}, mr={1324193}}
\bptok{imsref}%
\end{barticle}
%
\endbibitem

\bibitem[\protect\citeauthoryear{Hausman et~al.}{1991}]{Haus1991}
%
\begin{barticle}[mr]
\bauthor{\bsnm{Hausman},~\bfnm{Jerry~A.}\binits{J.~A.}},
\bauthor{\bsnm{Newey},~\bfnm{Whitney~K.}\binits{W.~K.}},
\bauthor{\bsnm{Ichimura},~\bfnm{Hidehiko}\binits{H.}} \AND
\bauthor{\bsnm{Powell},~\bfnm{James~L.}\binits{J.~L.}}
(\byear{1991}).
\btitle{Identification and estimation of polynomial errors-in-variables
models}.
\bjournal{J. Econometrics}
\bvolume{50}
\bpages{273--295}.
\bid{doi={10.1016/0304-4076(91)90022-6}, issn={0304-4076}, mr={1147115}}
\bptok{imsref}%
\end{barticle}
%
\endbibitem

\bibitem[\protect\citeauthoryear{Hu and Schennach}{2008}]{huschennachncme}
%
\begin{barticle}[mr]
\bauthor{\bsnm{Hu},~\bfnm{Yingyao}\binits{Y.}} \AND
\bauthor{\bsnm{Schennach},~\bfnm{Susanne~M.}\binits{S.~M.}}
(\byear{2008}).
\btitle{Instrumental variable treatment of nonclassical measurement error
models}.
\bjournal{Econometrica}
\bvolume{76}
\bpages{195--216}.
\bid{doi={10.1111/j.0012-9682.2008.00823.x}, issn={0012-9682}, mr={2374986}}
\bptok{imsref}%
\end{barticle}
%
\endbibitem

\bibitem[\protect\citeauthoryear{Huwang and Huang}{2000}]{huwangpoly}
%
\begin{barticle}[mr]
\bauthor{\bsnm{Huwang},~\bfnm{Longcheen}\binits{L.}} \AND
\bauthor{\bsnm{Huang},~\bfnm{Y.~H.~Steve}\binits{Y.~H.~S.}}
(\byear{2000}).
\btitle{On errors-in-variables in polynomial regression-{B}erkson case}.
\bjournal{Statist. Sinica}
\bvolume{10}
\bpages{923--936}.
\bid{issn={1017-0405}, mr={1787786}}
\bptok{imsref}%
\end{barticle}
%
\endbibitem

\bibitem[\protect\citeauthoryear{Hyslop and Imbens}{2001}]{hyslopncme}
%
\begin{barticle}[mr]
\bauthor{\bsnm{Hyslop},~\bfnm{Dean~R.}\binits{D.~R.}} \AND
\bauthor{\bsnm{Imbens},~\bfnm{Guido~W.}\binits{G.~W.}}
(\byear{2001}).
\btitle{Bias from classical and other forms of measurement error}.
\bjournal{J. Bus. Econom. Statist.}
\bvolume{19}
\bpages{475--481}.
\bid{doi={10.1198/07350010152596727}, issn={0735-0015}, mr={1963378}}
\bptok{imsref}%
\end{barticle}
%
\endbibitem

\bibitem[\protect\citeauthoryear{Lewbel}{1996}]{lewbel1996}
%
\begin{barticle}[author]
\bauthor{\bsnm{Lewbel},~\bfnm{A.}\binits{A.}}
(\byear{1996}).
\btitle{Demand estimation with expenditure measurement errors on the
left and
right hand side}.
\bjournal{Rev. Econom. Statist.}
\bvolume{78}
\bpages{718--725}.
\bptok{imsref}%
\end{barticle}
%
\endbibitem

\bibitem[\protect\citeauthoryear{Li}{2002}]{li1998}
%
\begin{barticle}[mr]
\bauthor{\bsnm{Li},~\bfnm{Tong}\binits{T.}}
(\byear{2002}).
\btitle{Robust and consistent estimation of nonlinear errors-in-variables
models}.
\bjournal{J. Econometrics}
\bvolume{110}
\bpages{1--26}.
\bid{doi={10.1016/S0304-4076(02)00120-3}, issn={0304-4076}, mr={1920960}}
\bptok{imsref}%
\end{barticle}
%
\endbibitem

\bibitem[\protect\citeauthoryear{Li and Vuong}{1998}]{livuong1998}
%
\begin{barticle}[mr]
\bauthor{\bsnm{Li},~\bfnm{Tong}\binits{T.}} \AND
\bauthor{\bsnm{Vuong},~\bfnm{Quang}\binits{Q.}}
(\byear{1998}).
\btitle{Nonparametric estimation of the measurement error model using multiple
indicators}.
\bjournal{J. Multivariate Anal.}
\bvolume{65}
\bpages{139--165}.
\bid{doi={10.1006/jmva.1998.1741}, issn={0047-259X}, mr={1625869}}
\bptok{imsref}%
\end{barticle}
%
\endbibitem

\bibitem[\protect\citeauthoryear{Mahajan}{2006}]{mahajanmisbin}
%
\begin{barticle}[mr]
\bauthor{\bsnm{Mahajan},~\bfnm{Aprajit}\binits{A.}}
(\byear{2006}).
\btitle{Identification and estimation of regression models with
misclassification}.
\bjournal{Econometrica}
\bvolume{74}
\bpages{631--665}.
\bid{doi={10.1111/j.1468-0262.2006.00677.x}, issn={0012-9682}, mr={2217611}}
\bptok{imsref}%
\end{barticle}
%
\endbibitem

\bibitem[\protect\citeauthoryear{Mallick, Hoffman and
Carroll}{2002}]{carrollnevada}
%
\begin{barticle}[mr]
\bauthor{\bsnm{Mallick},~\bfnm{Bani}\binits{B.}},
\bauthor{\bsnm{Hoffman},~\bfnm{F.~Owen}\binits{F.~O.}} \AND
\bauthor{\bsnm{Carroll},~\bfnm{Raymond~J.}\binits{R.~J.}}
(\byear{2002}).
\btitle{Semiparametric regression modeling with mixtures of {B}erkson and
classical error, with application to fallout from the {N}evada test site}.
\bjournal{Biometrics}
\bvolume{58}
\bpages{13--20}.
\bid{doi={10.1111/j.0006-341X.2002.00013.x}, issn={0006-341X}, mr={1891038}}
\bptok{imsref}%
\end{barticle}
%
\endbibitem

\bibitem[\protect\citeauthoryear{Nelder and Mead}{1965}]{nelderamoeba}
%
\begin{barticle}[author]
\bauthor{\bsnm{Nelder},~\bfnm{J.~A.}\binits{J.~A.}} \AND
\bauthor{\bsnm{Mead},~\bfnm{R.}\binits{R.}}
(\byear{1965}).
\btitle{A simplex method for function minimization}.
\bjournal{Computer Journal}
\bvolume{7}
\bpages{308--313}.
\bptok{imsref}%
\end{barticle}
%
\endbibitem

\bibitem[\protect\citeauthoryear{Newey}{1997}]{neweyseries}
%
\begin{barticle}[mr]
\bauthor{\bsnm{Newey},~\bfnm{Whitney~K.}\binits{W.~K.}}
(\byear{1997}).
\btitle{Convergence rates and asymptotic normality for series estimators}.
\bjournal{J. Econometrics}
\bvolume{79}
\bpages{147--168}.
\bid{doi={10.1016/S0304-4076(97)00011-0}, issn={0304-4076}, mr={1457700}}
\bptok{imsref}%
\end{barticle}
%
\endbibitem

\bibitem[\protect\citeauthoryear{Newey}{2001}]{Neweynliv}
%
\begin{barticle}[author]
\bauthor{\bsnm{Newey},~\bfnm{W.}\binits{W.}}
(\byear{2001}).
\btitle{Flexible simulated moment estimation of nonlinear errors-in-variables
models}.
\bjournal{Rev. Econom. Statist.}
\bvolume{83}
\bpages{616--627}.
\bptok{imsref}%
\end{barticle}
%
\endbibitem

\bibitem[\protect\citeauthoryear{Newey and Powell}{2003}]{neweynpiv}
%
\begin{barticle}[mr]
\bauthor{\bsnm{Newey},~\bfnm{Whitney~K.}\binits{W.~K.}} \AND
\bauthor{\bsnm{Powell},~\bfnm{James~L.}\binits{J.~L.}}
(\byear{2003}).
\btitle{Instrumental variable estimation of nonparametric models}.
\bjournal{Econometrica}
\bvolume{71}
\bpages{1565--1578}.
\bid{doi={10.1111/1468-0262.00459}, issn={0012-9682}, mr={2000257}}
\bptok{imsref}%
\end{barticle}
%
\endbibitem

\bibitem[\protect\citeauthoryear{Pope et~al.}{1995}]{popepartair}
%
\begin{barticle}[author]
\bauthor{\bsnm{Pope},~\bfnm{C.~A.}\binits{C.~A.}},
\bauthor{\bsnm{Thun},~\bfnm{M.~J.}\binits{M.~J.}},
\bauthor{\bsnm{Namboodiri},~\bfnm{M.~M.}\binits{M.~M.}} \betal{et~al.}
(\byear{1995}).
\btitle{Particulate air-pollution as a predictor of mortality in a
prospective-study of us adults}.
\bjournal{American Journal of Respiratory and Critical Care Medicine}
\bvolume{151}
\bpages{669--674}.
\bptok{imsref}%
\end{barticle}
%
\endbibitem

\bibitem[\protect\citeauthoryear{Samet et~al.}{2000}]{sametpartair}
%
\begin{barticle}[author]
\bauthor{\bsnm{Samet},~\bfnm{J.~M.}\binits{J.~M.}},
\bauthor{\bsnm{Dominici},~\bfnm{F.}\binits{F.}},
\bauthor{\bsnm{Curriero},~\bfnm{F.~C.}\binits{F.~C.}} \betal{et~al.}
(\byear{2000}).
\btitle{Fine particulate air pollution and mortality in 20 US cities,
1987--1994}.
\bjournal{New England Journal of Medicine}
\bvolume{343}
\bpages{1742--1749}.
\bptok{imsref}%
\end{barticle}
%
\endbibitem

\bibitem[\protect\citeauthoryear{Schennach}{2004}]{schennachnlme}
%
\begin{barticle}[mr]
\bauthor{\bsnm{Schennach},~\bfnm{Susanne~M.}\binits{S.~M.}}
(\byear{2004}).
\btitle{Estimation of nonlinear models with measurement error}.
\bjournal{Econometrica}
\bvolume{72}
\bpages{33--75}.
\bid{doi={10.1111/j.1468-0262.2004.00477.x}, issn={0012-9682}, mr={2031013}}
\bptok{imsref}%
\end{barticle}
%
\endbibitem

\bibitem[\protect\citeauthoryear{Schennach}{2007}]{schennachnlmeiv}
%
\begin{barticle}[mr]
\bauthor{\bsnm{Schennach},~\bfnm{Susanne~M.}\binits{S.~M.}}
(\byear{2007}).
\btitle{Instrumental variable estimation of nonlinear errors-in-variables
models}.
\bjournal{Econometrica}
\bvolume{75}
\bpages{201--239}.
\bid{doi={10.1111/j.1468-0262.2007.00736.x}, issn={0012-9682}, mr={2284741}}
\bptok{imsref}%
\end{barticle}
%
\endbibitem

\bibitem[\protect\citeauthoryear{Schennach}{2013}]{schennachberksupp}
%
\begin{bmisc}[author]
\bauthor{\bsnm{Schennach},~\bfnm{S.~M.}\binits{S.~M.}}
(\byear{2013}).
\btitle{Supplement to ``Regressions with Berkson errors in
covariates---A nonparametric approach.''
DOI:\doiurl{10.1214/13-AOS1122SUPP}}.
\bptok{imsref}%
\end{bmisc}
%
\endbibitem

\bibitem[\protect\citeauthoryear{Shen}{1997}]{Shensieve}
%
\begin{barticle}[mr]
\bauthor{\bsnm{Shen},~\bfnm{Xiaotong}\binits{X.}}
(\byear{1997}).
\btitle{On methods of sieves and penalization}.
\bjournal{Ann. Statist.}
\bvolume{25}
\bpages{2555--2591}.
\bid{doi={10.1214/aos/1030741085}, issn={0090-5364}, mr={1604416}}
\bptok{imsref}%
\end{barticle}
%
\endbibitem

\bibitem[\protect\citeauthoryear{Stram, Huberman and Wu}{2002}]{stramclasberk}
%
\begin{barticle}[pbm]
\bauthor{\bsnm{Stram},~\bfnm{Daniel~O.}\binits{D.~O.}},
\bauthor{\bsnm{Huberman},~\bfnm{Mark}\binits{M.}} \AND
\bauthor{\bsnm{Wu},~\bfnm{Anna~H.}\binits{A.~H.}}
(\byear{2002}).
\btitle{Is residual confounding a reasonable explanation for the apparent
protective effects of beta-carotene found in epidemiologic studies of lung
cancer in smokers?}
\bjournal{Am. J. Epidemiol.}
\bvolume{155}
\bpages{622--628}.
\bid{issn={0002-9262}, pmid={11914189}}
\bptok{imsref}%
\end{barticle}
%
\endbibitem

\bibitem[\protect\citeauthoryear{van~der Laan, Dudoit and
Keles}{2004}]{laanmlecv}
%
\begin{barticle}[mr]
\bauthor{\bparticle{van~der} \bsnm{Laan},~\bfnm{Mark~J.}\binits{M.~J.}},
\bauthor{\bsnm{Dudoit},~\bfnm{Sandrine}\binits{S.}} \AND
\bauthor{\bsnm{Keles},~\bfnm{Sunduz}\binits{S.}}
(\byear{2004}).
\btitle{Asymptotic optimality of likelihood-based cross-validation}.
\bjournal{Stat. Appl. Genet. Mol. Biol.}
\bvolume{3}
\bpages{Art. 4, 27 pp. (electronic)}.
\bid{doi={10.2202/1544-6115.1036}, issn={1544-6115}, mr={2101455}}
\bptok{imsref}%
\end{barticle}
%
\endbibitem

\bibitem[\protect\citeauthoryear{Wang}{2004}]{wangnlmbme}
%
\begin{barticle}[mr]
\bauthor{\bsnm{Wang},~\bfnm{Liqun}\binits{L.}}
(\byear{2004}).
\btitle{Estimation of nonlinear models with {B}erkson measurement errors}.
\bjournal{Ann. Statist.}
\bvolume{32}
\bpages{2559--2579}.
\bid{doi={10.1214/009053604000000670}, issn={0090-5364}, mr={2153995}}
\bptok{imsref}%
\end{barticle}
%
\endbibitem

\bibitem[\protect\citeauthoryear{Wang}{2007}]{wangunified}
%
\begin{barticle}[mr]
\bauthor{\bsnm{Wang},~\bfnm{Liqun}\binits{L.}}
(\byear{2007}).
\btitle{A unified approach to estimation of nonlinear mixed effects and
{B}erkson measurement error models}.
\bjournal{Canad. J. Statist.}
\bvolume{35}
\bpages{233--248}.
\bid{doi={10.1002/cjs.5550350203}, issn={0319-5724}, mr={2393607}}
\bptok{imsref}%
\end{barticle}
%
\endbibitem

\end{thebibliography}
\end{document}